\documentclass[08pt, a4paper, oneside]{article}
\usepackage[latin1]{inputenc}
\usepackage[english]{babel}
\usepackage{bm}
\usepackage{caption}
\usepackage{subcaption}
\usepackage{amsmath}
\usepackage{amsthm}
\usepackage{mathtools}
\usepackage{wasysym}
\usepackage{graphicx}
\usepackage{amstext}
\usepackage{color}
\usepackage{soul} 
\usepackage{sidecap}
\usepackage{amsbsy,amssymb}
\numberwithin{equation}{section}

\newcommand{\EE}{\ensuremath{\mathbf E}}
\newcommand{\FF}{\ensuremath{\mathbf F}}

\newcommand{\HH}{\ensuremath{\mathbf H}}
\newcommand{\II}{\ensuremath{\mathbf I}}

\newcommand{\MM}{\ensuremath{\mathbf M}}

\newcommand{\RR}{\ensuremath{\mathbf R}}

\newcommand{\WW}{\ensuremath{\mathbf W}}

\renewcommand{\aa}{\ensuremath{\mathbf a}}

\newcommand{\ii}{\ensuremath{\mathbf i}}

\newcommand{\mm}{\ensuremath{\mathbf m}}
\newcommand{\nn}{\ensuremath{\mathbf n}}

\newcommand{\uu}{\ensuremath{\mathbf u}}
\newcommand{\vv}{\ensuremath{\mathbf v}}

\newcommand{\zz}{\ensuremath{\mathbf z}}

\newcommand{\EEE}{\ensuremath{\mathbb E}}

\newcommand{\NNN}{\ensuremath{\mathbb N}}

\newcommand{\MMM}{\ensuremath{\mathbb M}}

\newcommand{\RRR}{\ensuremath{\mathbb R}}

\newcommand{\calA}{{\mathcal A}}
\newcommand{\calB}{{\mathcal B}}
\newcommand{\calC}{{\mathcal C}}

\newcommand{\calL}{{\mathcal L}}
\newcommand{\calX}{{\mathcal X}}
\newcommand{\calM}{{\mathcal M}}

\newcommand{\tr}{\mathrm{tr\,}}

\newtheorem{Theorem}{Theorem}

\newtheorem{Remark}{Remark}

\begin{document}

\title{Asymptotic analysis of hierarchical martensitic microstructure}

\author{Pierluigi Cesana\footnote{now at the Mathematical Institute, 
Woodstock Road, Oxford OX26GG, England}, Marcel Porta and Turab Lookman\\
\small Theoretical Division, Los Alamos National Laboratory, Los Alamos, NM 87545, USA}


 \maketitle

\begin{abstract}

We consider a hierarchical nested microstructure, which also contains a point of singularity (disclination) at the origin, observed in lead orthovanadate.  We  show how to  exactly compute the energy cost and associated displacement field within linearized elasticity by enforcing geometric compatibility of strains across interfaces of  the three-phase mixture of distortions (variants) in the microstructure.  We prove that the mechanical deformation is purely elastic and discuss the behavior of the system close to the origin.


\end{abstract}


\tableofcontents

\section{Introduction}

Solid-to-solid phase transformations are often accompanied by the formation of unusual and intriguing mixtures of phases at the mesoscale spanning  nanometers to microns in length scales \cite{Khachachuryan}. In the case of metallic alloys, below the transition temperature it is common to observe the coexistence  of fine layers of  martensitic or product twins (or "variants")  with the parent austenite phase of higher symmetry. Martensitic transformations are displacive, are driven by shear strain and/or shuffles  (intracell atomic displacements), and are invariably first-order in nature that leads to hysteresis and metastability \cite {BhattaBook}. The rich microstructure seen in high resolution electron microscopy (HREM) is often a manifestation of this metastability \cite{Manolikas, kitano}. Various approaches have been utilized over the last several decades to model the emergence of microstructure in martensites. They are largely variational using a free energy potential, and either use finite deformation, sharp interfaces and iteratively minimize a free energy, or start with random initial conditions and evolve a free energy potential  according to some dynamics \cite{Lookman01}. The traditional phase field approach is often used in the limit of small strains, and methods based on Ginzburg-Landau theory can be applied for small strains as well as finite deformation.

The study of mixtures in the framework of a variational setting traces back to the work of  Ball and James \cite{Ball87}. Their technique consists of matching different crystal phases, possibly at different scales, through geometric compatibility. In the literature, this idea has been widely applied in the analysis of periodic mixtures, the situation in which the physical and geometrical properties of these mixtures are repeated periodically in an elastic body.

Even though the case of periodic microstructure is of central importance in the modeling of composite materials, and it has become a classic subject of study, the situation related to more general microstructure, and specifically the case of self-similar fine hierarchies,  remains to be  fully explored.  In fact, there are a number of  fascinating open problems and questions,  among which is the study of a  large family of fine hierarchical structures observed in artificial polymers and biological materials (e.g., bones and leaves). Our emphasis will be on the family of heterogeneities in which there is  an interaction of topological singularities that leads to fascinating (non-periodic) microstructure.

In this paper we focus on the class of hexagonal-to-orthorhombic transformations where 
three equivalent stretching directions of the parent austenite phase  give rise to three orientations of the product martensite. Examples of compounds undergoing this transformation include the mineral  Mg-cordierite, Mg$_{2}$Al$_{4}$Si$_{5}$O$_{18}$, and Mg-Cd alloys \cite {kitano}. Closely related materials include those undergoing a hexagonal-to-monoclinic transformation, such as lead orthovanadate, Pb$_{3}$(VO$_{4}$)$_2$, and samarium sesquioxide, Sm$_{2}$O$_{3}$ \cite{Manolikas}.  
As the variants need to rotate to match at the domain walls, and the domain walls connecting the variants may intersect, these materials provide us with an excellent opportunity to study disclinations in crystals. Disclinations are formed when the nodes generated by the intersection of the domain walls do not close to an angle of 2$\pi$.  
 
\begin{figure}
        \centering
        \begin{subfigure}[b]{0.4\textwidth}
\includegraphics[width=6cm]{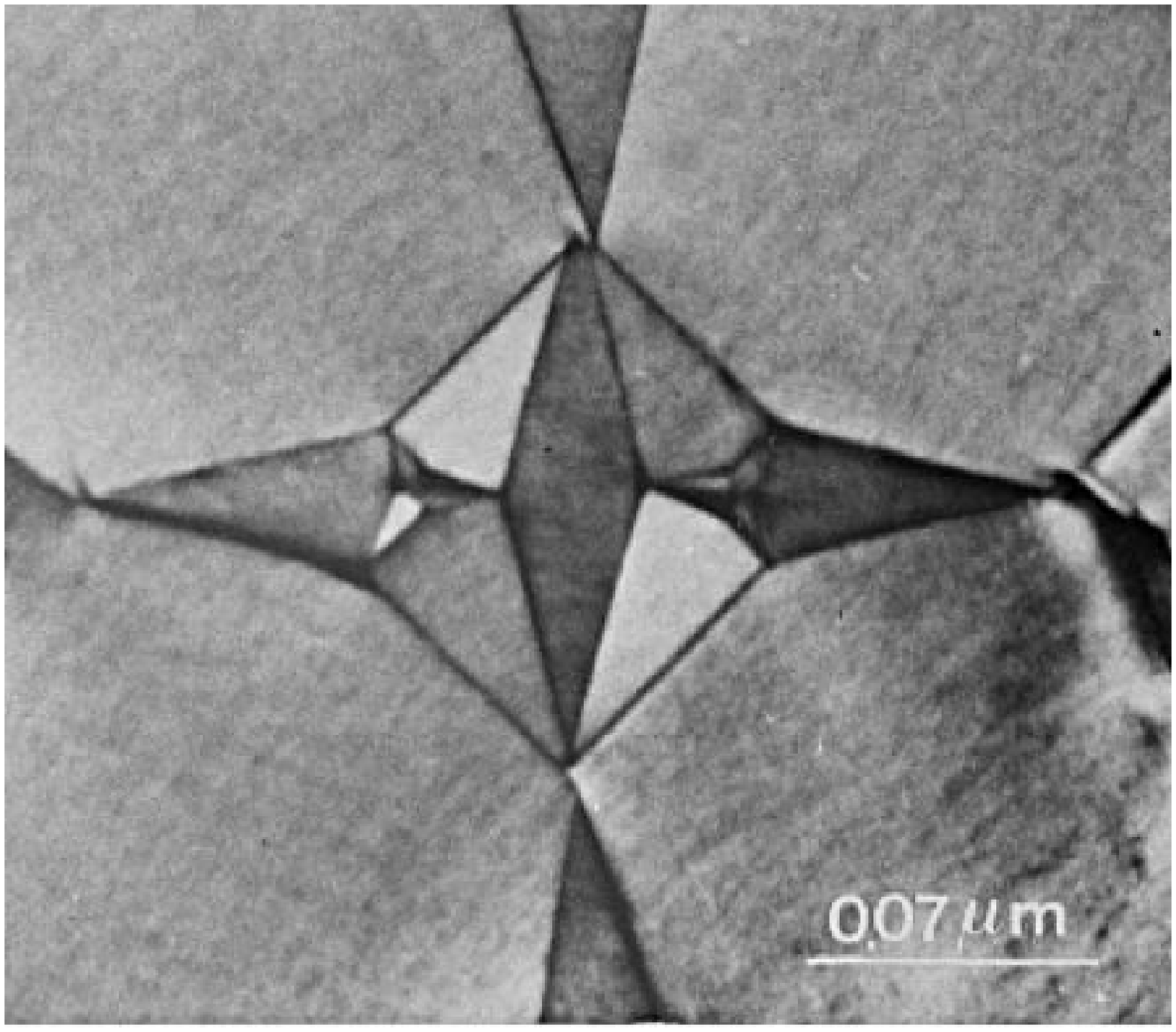}
                \caption{}
                \label{fig:a}
        \end{subfigure}%
        ~ 
        \begin{subfigure}[b]{0.3\textwidth}
\includegraphics[width=3.7cm]{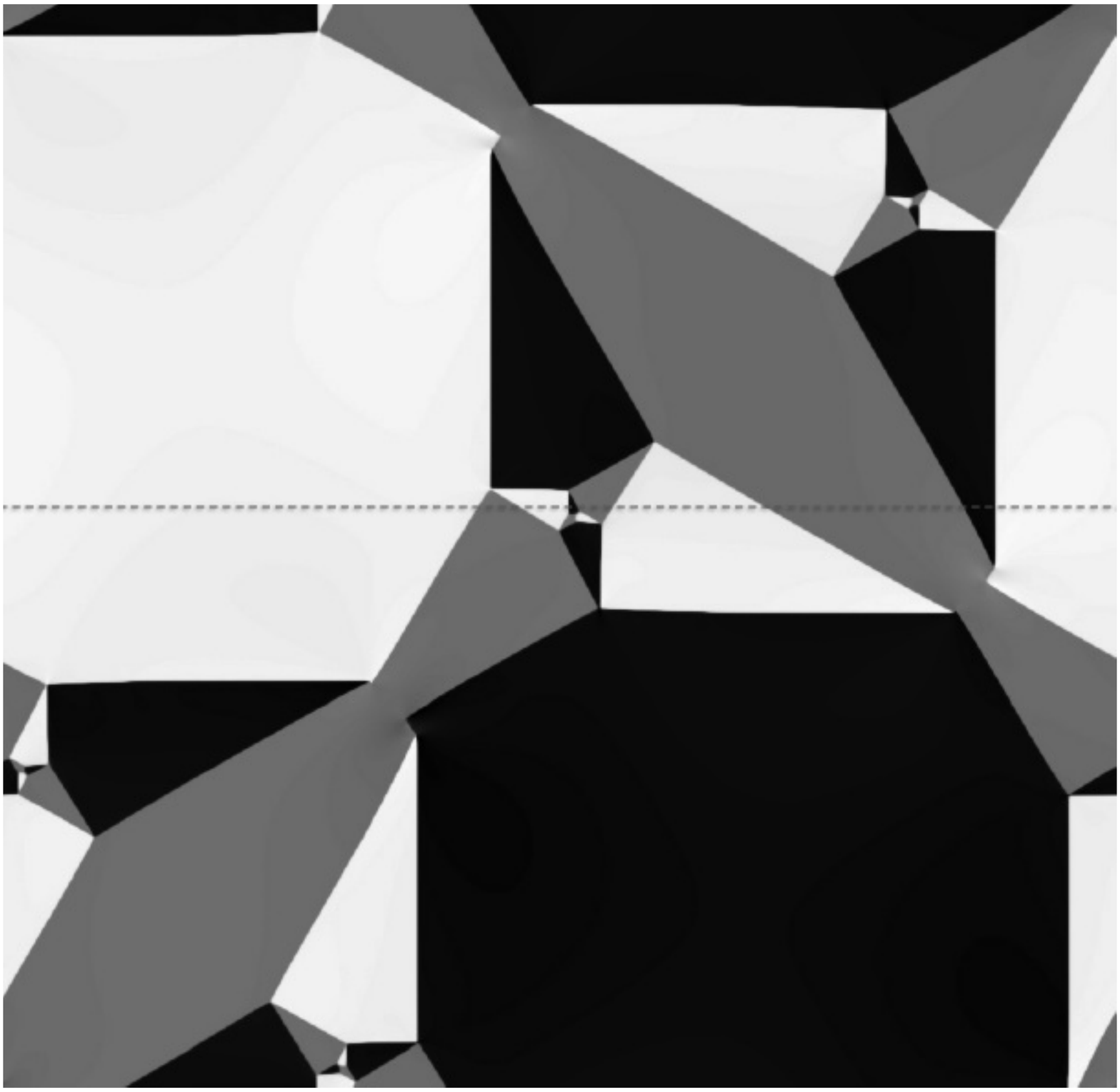}%
                \caption{}
                \label{fig:b}
        \end{subfigure}
        ~ 
\\
\quad
\\   
     \begin{subfigure}[b]{0.5\textwidth}
\centering
\includegraphics[width=7cm]{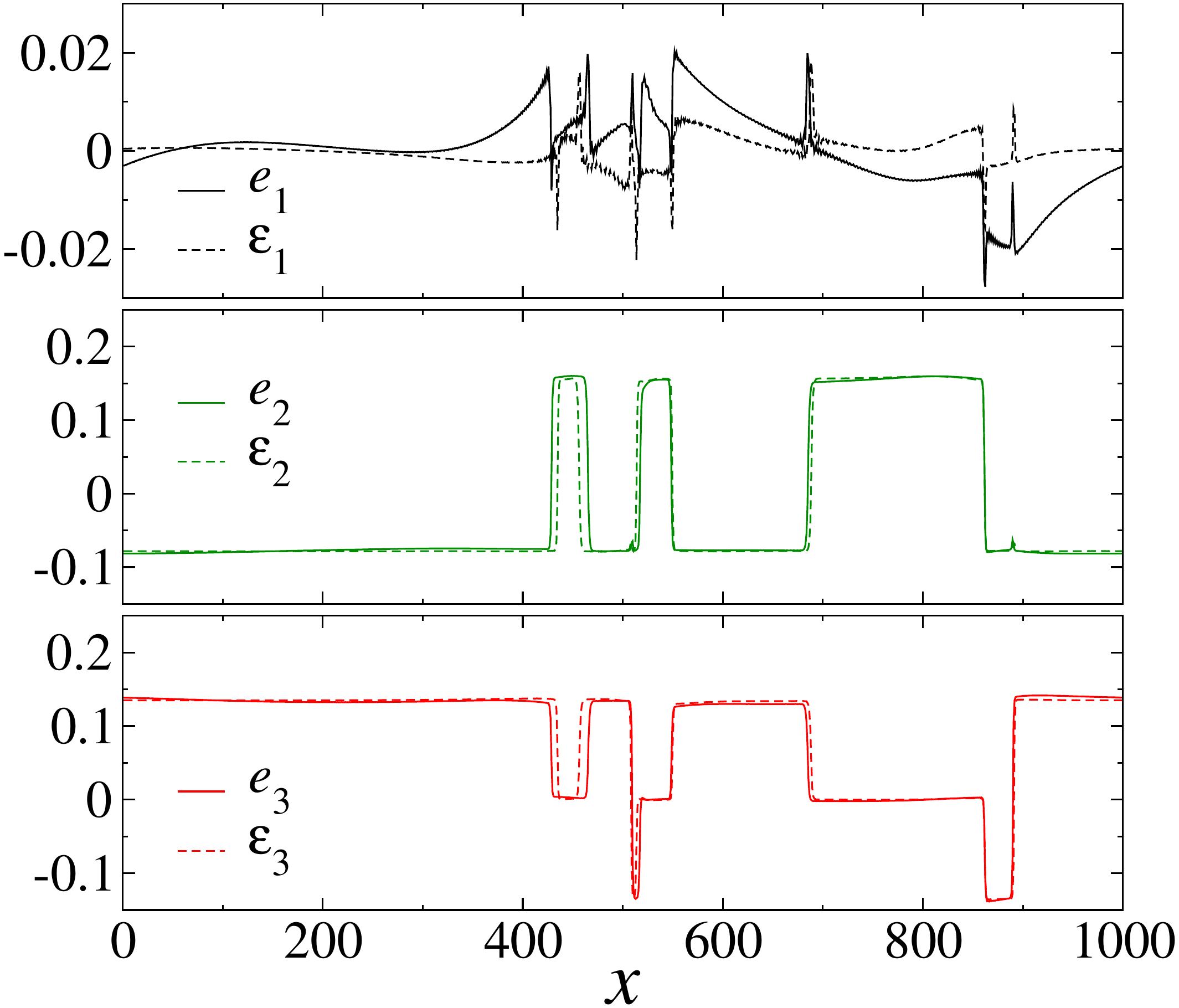}%
                \caption{}
                \label{fig:c}
        \end{subfigure}
\caption{(a): Experimental observation of the tripole-star microstructure in lead orthovanadate \cite {Manolikas}.
The two stars are modeled as rank-one three-phase martensitic mixtures. (b): Numerical solution of the mechanical equilibrium equations corresponding to the non-linear elastic model defined in Eqs. (\ref{13007012153}) and (\ref{1307011654}). (c): Comparison between the strain fields obtained numerically in the approximation of linear elasticity (dotted lines) and with the full nonlinear model (solid lines) (see also \cite{Lookman01}). The graphic shows the strain profiles along a line which crosses the tripole-star microstructure (indicated with a dotted line in (b)). In both cases the  choice of the parameters of the model yields $\epsilon\approx 0.156$.}
\label{1305141104}
\end{figure}

We will study the hexagonal-to-orthorhombic transformation in two dimensions (2D) for which the corresponding transformation is triangle-to-centered-rectangle. As the compounds consist of stacking or layering of tetrahedral units, the microstructure is essentially homogeneous perpendicular to the plane of the paper and therefore 2D is justified. One of the most intriguing microstructures observed in this transformation is the self-similarly nested tripole-star pattern (see Fig. \ref{1305141104}-(a)). These transformations have recently been the object of extensive numerical study. In Ref. \cite{Lookman01} the modeling is based on the minimization of a non-convex Ginzburg-Landau potential, both in the scenario of finite and infinitesimal elasticity.
It can be observed that in this microstructure the deformation gradient is (nearly) piecewise constant, with three possible strain values that correspond
to the three wells of the free energy density. As a first approximation, the sets where the deformation gradient is constant are particular $kyte$-shaped polygons. These polygons are all identical up to a rotation and a rescaling, close to the center of the star, and their measure tends to zero.

A thorough theory that treats microstructure and other phenomena, such as disclinations, dislocations, cavitations or cracks simultaneously, is to the best of our knowledge, still lacking.  Such a theory would provide an important tool to study various classes   of phenomena in materials science. Focusing on the case of linearized elasticity, the main result of this paper (Theorem \ref{THM}) is the exact computation of the zero-energy microstructure and of the displacement field that realizes it. 
In our construction we adapt the techniques of Ball and James \cite{Ball87} (see also \cite{Muller}) to the case of a non-periodic microstructure with a laborious construction from scratch. We note that the situation presented is a generalization of the case of simple laminates (matching of martensitic variants at the microscopic scale) or possibly, to the case of laminates-within-laminates. Indeed, the nesting of the microstructure requires a slightly more delicate construction consisting of simultaneous matching of phases across interfaces of polygons of different size at an infinite number of scales.
%
%
Furthermore, matching of geometrically compatible variants is typically achieved by piecewise affine displacement fields (laminates). Continuity of the displacement field is a key property of elastic models because it rules out the occurrence of irreversible phenomena such as cavitations and cracks.
In the present scenario, whether a continuous displacement can realize the microstructure depicted in Fig.\ref{1305141104}-(a) according to our model, is not a priori clear. The reason is that the center of the star may act as a topological singularity for the microstructure.
A by-product of our analysis in Theorem \ref{THM}, is that the displacement $\uu$ is indeed continuous and, consequently, it is a purely elastic displacement.

Whether adopting  a geometrically linearized model is physically sound for studying this type of microstructure is discussed in Ref. \cite{Lookman01} and in the final section of this paper. Although models in linearized elasticity have some intrinsic limitations that may rule out the understanding of irreversible or higher order phenomena, we remark here that the intrinsic simplicity of this model may shed some light on certain features that survive upon linearization, such as the geometry of the microstructure.
The analysis of the full non-linear model is left out from this article and is the object of ongoing research.

 The paper is organized as follows. In the next section we introduce the basic notation used through the paper and the basic concepts of finite elasticity. Moreover, we
present the mechanical model used in Ref. \cite{Lookman01} and its geometrical linearization obtained in the asymptotic expansion of small displacements. In Section \ref{tripole}  we present the full analysis of the geometry of the problem and the construction of a possible displacement field $\uu$ which reproduces a microstructure similar to the one observed in Fig. \ref{1305141104}-(a) (Theorem \ref{THM}). Finally, Section \ref{1307112003} is devoted to the physical interpretation of the results contained in Theorem \ref{THM} and to the presentation of some open problems.


%
%
%
%
%
%
%
%
%

\section{Background}

\subsection{Notation}

We gather here the main symbols and the notation used throughout
the paper. Our main references are \cite{BhattaBook} and \cite{Ciarlet}. Let $\NNN$ and $\RRR$ denote the set of natural and
real numbers respectively.
For any integer $n$, $\RRR^{n}$ is the
space of $n$-dimensional vectors with canonical basis $\{\ii_i\},
i=1,\dots n$ with origin $O=(0,\dots,0)$ and $\MMM^{n\times n}$
the space of square real matrices. The determinant, the trace and
the transpose of the matrix $\FF$ in $\MMM^{n\times n}$ are
denoted by $\det\FF$, $\tr\FF$, $\FF^T$ respectively. 
We endow
$\MMM^{n\times n}_{}$ with the usual inner product
$\FF:\MM:=\tr(\FF\MM^T)=\sum_{ij}F_{ij}M_{ij}$ and the corresponding
norm $|\MM|:=(\MM:\MM)^{1/2}$.
Here
$M_{ij},F_{ij}$ are the cartesian components of $\MM$ and $\FF$.
The identity in $\MMM^{n\times n}$ is denoted with $\II$ with components $\delta_{ij}$. 
We have the orthogonal decomposition of a matrix $\FF\in\MMM^{n\times n}$:
$\MM=\MM_{sym}+\MM_{skew}$
where $\MM_{sym}:=(\MM+\MM^T)/2$ and $\MM_{skew}:=(\MM-\MM^T)/2$. By further decomposition of $\MM_{sym}$ into its deviatoric and spherical part we have 
\begin{eqnarray}\label{1306301820}
\MM=\MM_{dev}+\MM_{sph}+\MM_{skew}
\end{eqnarray}
where $\MM_{sph}:=(\frac{\tr\MM}{n})\II$ and
$\MM_{dev}:=\MM_{sym}-\MM_{sph}$.
 
We denote with $B_{R}(0,0)$ the disk in $\RRR^2$ with center in the origin and with radius equal to $R$:
\begin{eqnarray}
B_{R}(0,0):=\{(x,y)\in\RRR^2:x^2+y^2<R^2\}.
\end{eqnarray}
Let $\Omega$ be an open and bounded subset of $\RRR^n$. We denote with $C(\Omega,\RRR^n)$ the space of continuous vector fields and with $C(\overline{\Omega},\RRR^n)$ those maps in $C(\Omega,\RRR^n)$ which are continuous on the closure of $\Omega$.
Letting $p\in[1,+\infty)$, we introduce $L^{p}(\Omega)$,
the space of measurable functions $u:\Omega\mapsto\RRR$ such that
$\int_{\Omega}|u|^pdx<+\infty$. Analogously,
$L^{p}(\Omega,\RRR^n)$ and $L^{p}(\Omega,\MMM^{n\times n})$,
respectively the spaces of vectors or matrices with components in
$L^p(\Omega)$. Then,
$W^{1,p}(\Omega,\RRR^n)$ is
the spaces of vector-valued $L^p$-functions
whose gradient has $L^p$-integrable components.  
The space $W^{1,\infty}(\Omega,\RRR^n)$ is
that of vectors with essentially bounded components whose gradient has essentially bounded components. For $\Omega$ regular enough we identify $W^{1,\infty}(\Omega,\RRR^n)$ with the space of Lipschitz functions over $\Omega$.
Other spaces of functions may be defined when encountered throughout the paper.

\subsection{Finite elasticity}

According to Ref. \cite{BhattaBook} Chapter 2 we introduce the deformation gradient $\FF$ whose components are defined as
$$
F_{ij}=\frac{\partial \calX_{i}}{\partial x_{j}},
$$
where $\calX_i$ is the $i$th component of the position vector of a mass element  in the current configuration and $x_j$ is the $j$th component of its position vector in the initial reference configuration. In the case $n=2$ we let $x_1\equiv x$, $x_2\equiv y$.
The deformation gradient can be written in terms of the displacement gradient as
\begin{eqnarray}\label{1305131210}
F_{ij}=\frac{\partial\bigl((\calX_{i}- x_i)+x_i\bigr)}{\partial x_{j}}=\frac{\partial (u_{i}+ x_i)}{\partial x_{j}}=\frac{\partial u_i}{\partial x_j}+\delta_{ij},
\end{eqnarray}
where $u_i$ are the components of the displacement vector $\uu$.
The Lagrangian strain tensor is defined as
$\EEE=\frac{1}{2}(\FF^T\FF-\II)$.
In components, the Lagrangian strain tensor reads 
\begin{eqnarray}\label{}
\EEE_{ij}=\frac{1}{2}\Bigl(\frac{\partial u_i}{\partial x_j}+\frac{\partial u_j}{\partial x_i}+\sum_k\frac{\partial u_k}{\partial x_i}\frac{\partial u_k}{\partial x_j}\Bigr).
\end{eqnarray}
%
%
Finally, we present the symmetry-adapted Lagrangian strains expressed, as in \cite[Section 6.1]{Lookman01}, as functions of the displacement gradient:
%
%
%
%
\begin{gather}
e_1:=\frac{\EEE_{xx}+\EEE_{yy}}{2}=\frac{1}{4}(F_{11}^2+F_{21}^2+F_{12}^2+F_{22}^2-2)=\nonumber\\
\frac{1}{2}\Bigl[\frac{\partial u_x}{\partial x}+\frac{\partial u_y}{\partial y}+\frac{1}{2}\Bigl(\frac{\partial u_x}{\partial x}\Bigr)^2+\frac{1}{2}\Bigl(\frac{\partial u_y}{\partial y}\Bigr)^2+\frac{1}{2}\Bigl(\frac{\partial u_x}{\partial y}\Bigr)^2+\frac{1}{2}\Bigl(\frac{\partial u_y}{\partial x}\Bigr)^2\Bigr],\nonumber\\
e_2:=\frac{\EEE_{xx}-\EEE_{yy}}{2}=\frac{1}{4}(F_{11}^2+F_{21}^2-F_{12}^2-F_{22}^2)=\nonumber\\
\frac{1}{2}\Bigl[\frac{\partial u_x}{\partial x}-\frac{\partial u_y}{\partial y}+\frac{1}{2}\Bigl(\frac{\partial u_x}{\partial x}\Bigr)^2-\frac{1}{2}\Bigl(\frac{\partial u_y}{\partial y}\Bigr)^2-\frac{1}{2}\Bigl(\frac{\partial u_x}{\partial y}\Bigr)^2+\frac{1}{2}\Bigl(\frac{\partial u_y}{\partial x}\Bigr)^2\Bigr],\nonumber\\
e_3:=\frac{\EEE_{xy}+\EEE_{yx}}{2}=\frac{1}{2}(F_{11}F_{12}+F_{21}F_{22})=\frac{1}{2}\Bigl[\frac{\partial u_x}{\partial y}+\frac{\partial u_y}{\partial x}+\frac{\partial u_x}{\partial x} \frac{\partial u_x}{\partial y} + \frac{\partial u_y}{\partial x} \frac{\partial u_y}{\partial y}  \Bigr].
\label{1305131522}\end{gather}
%
%

\subsection{The Ginzburg-Landau model for the triangle-to-centered-rectangle transformation}
\label{1307071856}

In this Section we briefly review the Ginzburg-Landau model for the 2D triangle-to-centered-rectangle (TR) transformation introduced in Refs.
\cite{Lookman01,Jacobs,Lookman03}. The TR transformation is the 2D version of the 3D hexagonal-to-orthorhombic transformation that occurs in materials such as the MgCd ordered alloy. A peculiar aspect of these materials is that they are rare examples of disclinations in crystals. In the 2D setting, disclinations are point defects that result when 
multiple twin boundaries intersect. The matching condition between two martensite  variants requires the crystal lattice to rotate, and the total rotation angle along a closed circuit that contains the intersection of the twin boundaries (disclination) is, in general, nonzero. This  requires a  stretch of the lattice, additional to the transformation strain, in order to maintain its coherency.
The disclination is characterized by the disclination angle, which is the total rotation of the crystal lattice along the closed circuit that contains the defect.

A special case of disclination and the most interesting pattern in the microstructure of the MgCd alloy is the self-similarly nested tripole-star disclination,
also observed in lead orthovanadate, $\textrm{Pb}_3(\textrm{VO}_4)_2$, undergoing a trigonal-to-monoclinic transformation \cite{Manolikas}. 
This pattern was analyzed in detail by Kitano and Kifune \cite{kitano}.  
In this microstructure the intersection of twin boundaries do not require a stretch of the lattice for their matching, and therefore they are not disclinations. However, the star pattern is itself a disclination, as the total rotation angle along a closed circuit containing the center of the star is nonzero.

In order to study this microstructure, 
in this paper we adopt the Ginzburg-Landau approach appropriate for the TR transformation.
The corresponding Landau free energy density, $\psi_{L}$, is a function of  the symmetry adapted components of the Lagrangian strain tensor, and has the functional form
\begin{figure}[h!]
\centering%
\includegraphics[width=6cm]{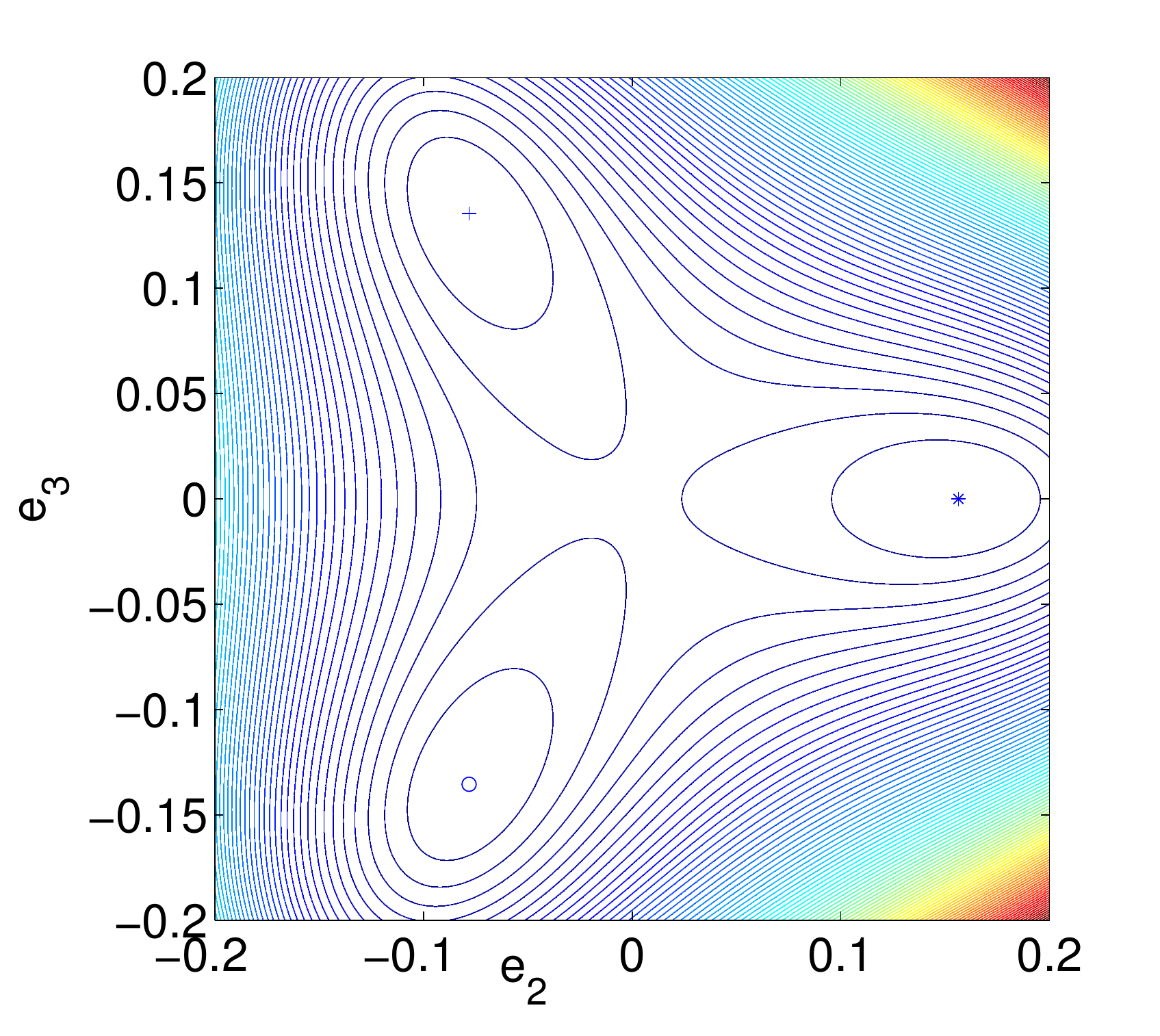}
\caption{Level curves of $\psi_{L}$ for $e_1\equiv 0$. Here $\calA=1$, $\calB=-30$, $\calC=200$, $T=0.8T_C$, $T_C=1$ ($\epsilon\approx 0.156$).}\label{quadratitris}
\end{figure}
\begin{gather}\label{13007012153}
\psi_{L}(e_1,e_2,e_3)=\frac{\calA_1}{2} e_1^2+\frac{\calA}{2}(T-T_C)(e_2^2+e_3^2)+\frac{\calB}{3}(e_2^3-3e_2 e_3^2)+\frac{\calC}{4}(e_2^2+e_3^2)^2+\calM.
\end{gather}
The phase transition is first order for ${\cal B} \ne 0$ and the transition temperature is $T_0 = T_C + 2\calB^2/9\calA\calC$.
For $T<T_0$ the free energy density has three absolute minima, corresponding to the three martensitic variants, with strain values,
\begin{gather}\label{1305272015}
\left.
\begin{array}{cccc}
 e_1=0,\qquad  e_2=\epsilon,\qquad & e_3=0\\
 e_1=0,\qquad  e_2=-\frac{\epsilon}{2},\qquad & e_3=\sqrt{3}e_2\\
 e_1=0,\qquad  e_2=-\frac{\epsilon}{2},\qquad & e_3=-\sqrt{3}e_2,\\
\end{array}
\right.
\end{gather}
and 
$$
\epsilon=\frac{-\calB+\sqrt{\calB^2-4\calC\calA(T-T_C)}}{2\calC}.
$$
For ${\cal B}<0$, $\epsilon$ is a positive (real) number. 
The parameter $\calM$ in (\ref{13007012153}) is chosen so that  $\min\psi_{L}=0$. Thus, 
$$
\calM=-\Bigl[\frac{\calA}{2}(T-T_C)\epsilon^2+\frac{\calB}{3}\epsilon^3+\frac{\calC}{4}\epsilon^4\Bigr].
$$
Moreover, variations in the order parameter strain fields are incorporated through a gradient term.
Thus, the total free energy density considered is
\begin{equation}
\psi_{GL}=\psi_{L}+ \frac{K}{2}\bigl(|\nabla e_2(x,y)|^2+|\nabla e_3(x,y)|^2\bigr).
\label{fGL}
\end{equation}
The approach adopted in Ref. \cite{Lookman01} consists of minimizing the total energy of the system,
\begin{eqnarray}\label{1307011654}
\Psi= \int_{\Omega}\Bigl[\psi_{L}\bigl(e_1(x,y),e_2(x,y),e_3(x,y)\bigr)+ \frac{K}{2}\bigl(|\nabla e_2(x,y)|^2+|\nabla e_3(x,y)|^2\bigr) \Bigr]dxdy,
\end{eqnarray}
with the constraint of compatibility of the strain fields. This can be done using a relaxational algorithm or solving the Euler-Lagrange equations associated with (\ref{1307011654}). 

In Fig.\ref{1305141104}-(b) we show a snapshot of a microstructure with several tripole-star patterns obtained as a solution of  the mechanical equilibrium equations using an iterative spectral method \cite{Lookman01}. The microstructure is obtained with $\epsilon\approx 0.156$ and only the strain component $e_3$ is shown, with light (dark) greyscale corresponding to positive (negative) values of $e_3$.
Notice that $e_3$ is (approximately) constant over each subregion of the same color. The same behavior is also observed for $e_2$. On the contrary, the strain component $e_1$ is (approximately) zero within the subregions, whereas nonzero at the boundaries of the subregions (twin boundaries).
While the effect of the regularizing term $\frac{K}{2}(|\nabla e_2|^2+|\nabla e_3|^2)$ is to smooth the twin boundaries separating the different subregions, the values of the symmetric Lagrangian strain components within each subregion are determined by 
$\psi_{L}$.
 
The full non-linear model (Eq. (\ref{1307011654}))  is
appropriate for the analysis of disclinations in crystals, as it can distinguish the different types of interfaces present in the microstructure \cite{curnoe} which are the origin of these defects.
Although a  geometric linearization of this model does not capture the existence of the disclination, the linear model is still able to reproduce
the tripole-star pattern \cite{Lookman01}.
In the next Subsection we recall the approximations made in geometric linear elasticity and apply them to the Ginzburg-Landau model introduced in this Section.

\subsection{Geometric linearization of the Ginzburg-Landau model}
\label{1211121518}

In the hypothesis of small displacements (and small displacement gradients) one can actually neglect the higher order terms in (\ref{1305131522}) yielding
\begin{eqnarray}\label{1307081523}
e_1=\frac{\EEE_{xx}+\EEE_{yy}}{2}\approx\frac{1}{2}\Bigl[\frac{\partial u_x}{\partial x}+\frac{\partial u_y}{\partial y}\Bigr]:=\varepsilon_1(x,y),\nonumber\\
e_2=\frac{\EEE_{xx}-\EEE_{yy}}{2}\approx\frac{1}{2}\Bigl[\frac{\partial u_x}{\partial x}-\frac{\partial u_y}{\partial y}\Bigr]:=\varepsilon_2(x,y),\nonumber\\
e_3=\frac{\EEE_{xy}+\EEE_{yx}}{2}\approx\frac{1}{2}\Bigl[\frac{\partial u_x}{\partial y}+\frac{\partial u_y}{\partial x}\Bigr]:=\varepsilon_3(x,y).
\end{eqnarray} 
Upon introduction of $\EE=\EE(\nabla\uu):=\nabla\uu_{sym}$ the (linearized) mechanical strain tensor, we have that $\varepsilon_2$ and $\varepsilon_3$ are the components of the deviatoric part of $\EE$ and $\varepsilon_1$ is the component of the spherical part of $\EE$:
\begin{eqnarray}\label{1307041431}
\EE=\frac{1}{2}(\nabla\uu+\nabla\uu^T)=\left(
\begin{array}{ccc}
\varepsilon_{2} & \varepsilon_3  \\
\varepsilon_3 & -\varepsilon_{2}
\end{array} \right)+\varepsilon_1\II.
\label{1110120953}\end{eqnarray}
In the linearized elasticity scenario, the new mechanical model is simply obtained by plugging $\varepsilon_1, \varepsilon_2$ and $\varepsilon_3$ instead of $e_1,e_2$ and $e_3$ respectively into the Ginzburg-Landau free energy density defined in Eqs. (\ref{13007012153}) and (\ref{fGL}). The Landau free energy density, which we still denote by $\psi_{L}$, reads
\begin{gather}\label{1307012156}
\psi_{L}(\varepsilon_1,\varepsilon_2,\varepsilon_3)=\frac{\calA_1}{2} \varepsilon_1^2+\frac{\calA}{2}(T-T_C)(\varepsilon_2^2+\varepsilon_3^2)+\frac{\calB}{3}(\varepsilon_2^3-3\varepsilon_2 \varepsilon_3^2)+\frac{\calC}{4}(\varepsilon_2^2+\varepsilon_3^2)^2+\calM.
\end{gather}
Minimization of (\ref{1307012156}) yields, in terms of strain matrices written as in (\ref{1307041431}), the three minimum points
\begin{eqnarray}\label{1305131529}
\EE_1= \epsilon\left(
\begin{array}{ccc}
1 & 0  \\
0 & -1
\end{array} \right),
\quad
\EE_2= \epsilon\left(
\begin{array}{ccc}
-\frac{1}{2} & \frac{\sqrt{3}}{2}  \\
\frac{\sqrt{3}}{2} & \frac{1}{2}
\end{array} \right),
\quad
\EE_3= \epsilon\left(
\begin{array}{ccc}
-\frac{1}{2} & -\frac{\sqrt{3}}{2}  \\
-\frac{\sqrt{3}}{2} & \frac{1}{2}
\end{array} \right).
\label{1110120953}
\end{eqnarray}
Notice that even in the linearized model $\psi_{L}$ has three distinct minimum points. 
Importantly: 
\begin{eqnarray}\label{1306261933}
|\EE_1|^2=|\EE_2|^2=|\EE_3|^2=2\epsilon^2.
\end{eqnarray}

\section{The tripole-star pattern}\label{tripole}

Minimization  of the total energy 
\begin{eqnarray}\label{1307041438}
\int_{\Omega_R}\Bigl[\psi_{L}\bigl(\varepsilon_1(x,y),\varepsilon_2(x,y),\varepsilon_3(x,y)\bigr) + \frac{K}{2}\bigl(|\nabla \varepsilon_2(x,y)|^2+|\nabla \varepsilon_3(x,y)|^2\bigr)\Bigr]dxdy
\end{eqnarray}
yields numerical solutions for the strain fields which are, at least from a qualitative point of view, analogous to those obtained for the full nonlinear model {\cite{Lookman01}}.
Briefly, the term $\psi_{L}$ enforces the symmetrized gradient of $\uu$ to be very close to the either $\EE_1,\EE_2$ or $\EE_3$ while the presence of $\frac{K}{2}\bigl(|\nabla \varepsilon_2|^2+|\nabla \varepsilon_3|^2\bigr)$ prevents the formation of steep interfaces between the three phases. 
Additional models are available in the literature which 
penalize the length of the interfaces and at the same time
allow jump discontinuities of the variants. 
%
%
We refer specifically to a forthcoming paper of Dr. A. Ruland for an analysis of a model with interfacial energy terms for the tripole-star considered in this paper.
In Fig.\ref{1305141104}-(c) we compare the strain fields obtained with the linearized model with the strain fields obtained with the full nonlinear model.
%
%
As previously noted and pointed out in Ref. \cite{Lookman01}, a Ginzburg-Landau type of model can offer an excellent insight even though being defined for infinitesimal displacements.
Inspired by these observations, we analyze in detail the linearized model.
We present the construction of a microstructure which is very similar to both the experimental observation and to the numerical solution plotted in Fig. \ref{1305141104} for the full-elasticity model by allowing sharp interfaces, that is, when we set $K\equiv 0$ in (\ref{1307041438}).

For the reader's convenience, our construction is split in two parts. In the next subsection, we present the geometry of the microstructure. 
Subsequently, we describe the requirements that a map $\uu$ has to fulfill in order to produce a microstructure which is consistent with the experimental observation (see Fig.  \ref{1305141104}-(b)). Theorem \ref{THM}, which is the core of this paper, contains the detailed construction of $\uu$.

\subsection{Geometry of the microstructure}\label{1003151615}

\begin{figure}[h!]
\centering
\includegraphics[trim = 60mm 00mm 60mm 00mm, clip=true, width=6.0cm]{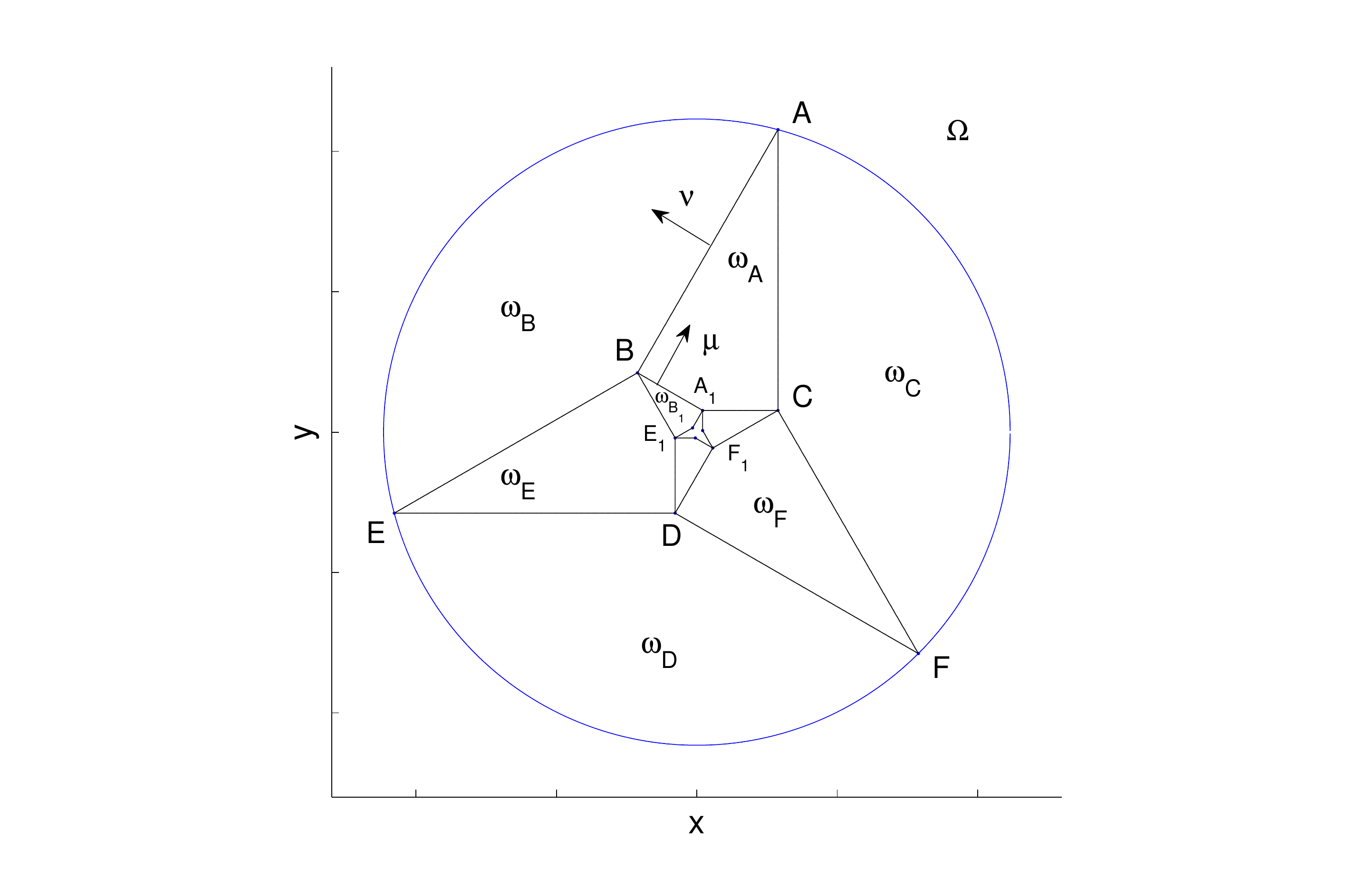}
\quad
\includegraphics[trim = 20mm 00mm 20mm 00mm, clip=true, width=6.0cm]{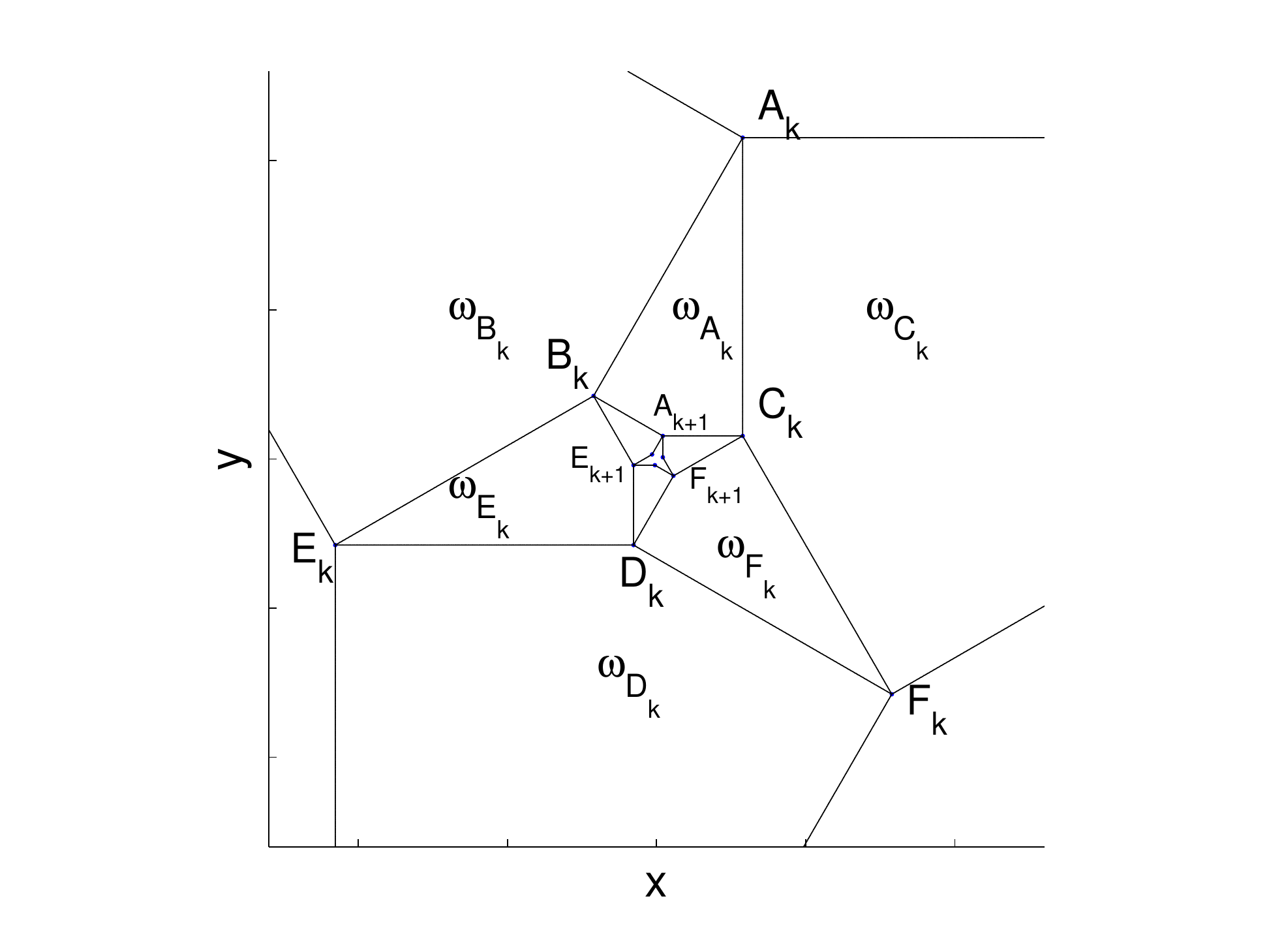}
\caption{LEFT: $\Omega$ and $\omega_0$. RIGHT: close-up on $\omega_{A_k}$, $\omega_{E_k}$ and $\omega_{F_k}$, $k\geq 1$. Notice that the points $B_{k-1}$, $C_{k-1}$ and $D_{k-1}$ are not shown in the picture.}\label{1111121420}
\end{figure}

\begin{figure}[h!]
\centering
\includegraphics[trim = 60mm 00mm 60mm 00mm, clip=true, width=8.5cm]{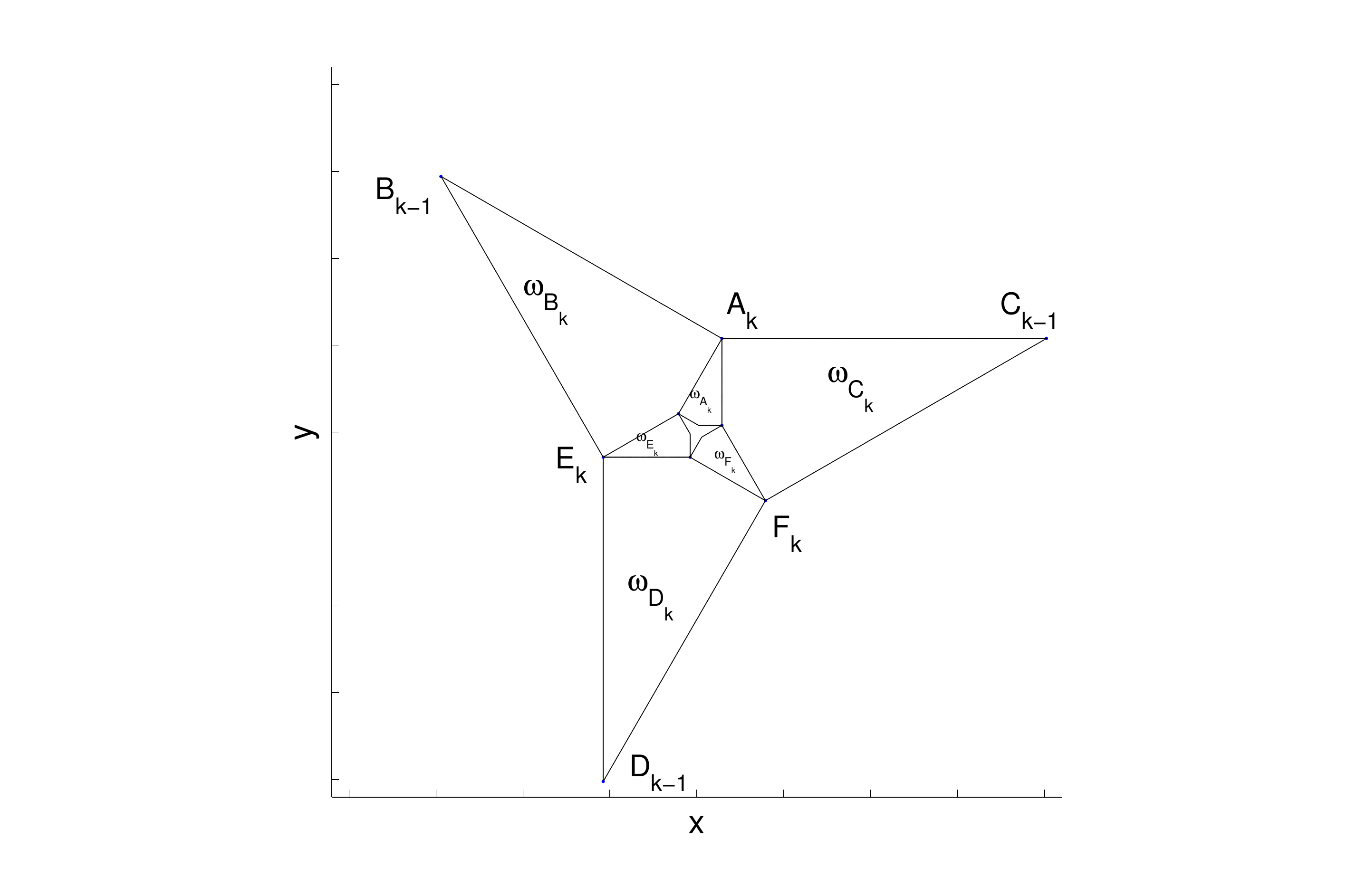}
\caption{$\omega_k$, $k\geq 1$.}\label{1307081558}
\end{figure}

We assume as our reference configuration the set $\Omega:=B_{R}(0,0)$, the disk with center in the origin and of radius equal to $R$.
We define $R:=L\sqrt{\frac{2}{3}+\frac{1}{\sqrt{3}}}$ where $L$ is a positive constant. In turn, $L$ can be regarded as the characteristic length of the system.
To reproduce the geometry of Fig. \ref{1305141104} we combine a family \textit{kyte}-shaped rhomboids that cover $\Omega$  (see Fig. \ref{1111121420}-LEFT).
Let $t:=\tan(\pi/12)$ and $k\in\{0\}\cup\NNN$.
We define the vertices of the rhomboids (here and in what follows, when $k=0$ we identify $A_0\equiv A$ and similarly for $B,C,D,E,F$):
\begin{eqnarray}\label{1307151522}
\left\{
\begin{array}{cccccc}
A &= &(x_A,y_A) &= &L(\frac{1}{2\sqrt{3}},\frac{1}{2}+\frac{1}{\sqrt{3}})  \\
B &= &(x_B,y_B) &= &L(\frac{1-\sqrt{3}}{2\sqrt{3}},-\frac{1-\sqrt{3}}{2\sqrt{3}})  \\
C &= &(x_C,y_C) &= &L(\frac{\sqrt{3}}{6},\frac{1}{\sqrt{3}}-\frac{1}{2})  \\
D &= &(x_D,y_D) &= &L(\frac{1}{2}-\frac{1}{\sqrt{3}},-\frac{1}{2\sqrt{3}})  \\
E &= &(x_E,y_E) &= &L( \frac{-2-\sqrt{3}}{2\sqrt{3}},-\frac{1}{2\sqrt{3}})  \\
F &= &(x_F,y_F) &= &L(\frac{1+\sqrt{3}}{2\sqrt{3}},-\frac{1+\sqrt{3}}{2\sqrt{3}})  \\
\end{array} \right.
\qquad
\left\{
\begin{array}{cccccc}
A_k &= &(x_{A_k},y_{A_k}) &= &t^{2k}A \\
B_k &= &(x_{B_k},y_{B_k}) &= &t^{2k}B \\
C_k &= &(x_{C_k},y_{C_k}) &= &t^{2k}C \\
D_k &= &(x_{D_k},y_{D_k}) &= &t^{2k}D \\
E_k&= &(x_{E_k},y_{E_k}) &= &t^{2k}E\\
F_k &= &(x_{F_k},y_{F_k}) &= &t^{2k}F\\
\end{array} \right.
\end{eqnarray}
In table $\ref{1211091745}$ we report the complete description of the segments that define the sides of the polygons.
\begin{table}[h]
\centering
\begin{tabular}{|c|l|r|}\hline
 $\quad$ Interface $\quad$  & Equation of the interface  \\ \hline
$\overline{B_kA_k}$  &$y=g_k^{BA}(x):=\sqrt{3}x+\frac{L}{\sqrt{3}}t^{2k},\qquad\qquad\qquad\,  x_{B_k}<x<x_{A_k}$\\ \hline
$\overline{C_kA_k}$  &$x=g_k^{CA}(y):=\frac{L}{2\sqrt{3}}t^{2k},\qquad\qquad\qquad\qquad\quad y_{C_k}<y<y_{A_k}$\\\hline
$\overline{C_kF_k}$  &$y=g_k^{CF}(x):=-\sqrt{3}x+\frac{L}{\sqrt{3}}t^{2k},\qquad\qquad\quad\,\, x_{C_k}<x<x_{F_k}$\\ \hline
$\overline{D_kF_k}$  &$y=g_k^{DF}(x):=-\frac{1}{\sqrt{3}}x-\frac{L}{3}t^{2k},\qquad\qquad\quad\,\,\,\,\, x_{D_k}<x<x_{F_k}$\\\hline
$\overline{E_kD_k}$  &$y=g_k^{ED}(x):=-\frac{L}{2\sqrt{3}}t^{2k},\qquad\qquad\qquad\qquad x_{E_k}<x<x_{D_k}$\\\hline
$\overline{E_kB_k}$  &$y=g_k^{EB}(x):=\frac{1}{\sqrt{3}}x+\frac{L}{3}t^{2k},\qquad\qquad\qquad\,\,\, x_{E_k}<x<x_{B_k}$\\\hline
$\overline{B_kA_{k+1}}$  &$y=h^{BA}_k(x):=-\frac{1}{\sqrt{3}}x+(\frac{2}{3}-\frac{1}{\sqrt{3}})Lt^{2k},\qquad x_{B_k}<x<x_{A_{k+1}}$\\ \hline
$\overline{A_{k+1}C_k}$  &$y=h^{AC}_k(x):=(\frac{1}{\sqrt{3}}-\frac{1}{2})Lt^{2k},\qquad\qquad\quad x_{A_{k+1}}<x<x_{C_k}$\\\hline
$\overline{F_{k+1}C_k}$  &$y=h^{FC}_k(x):=\frac{x}{\sqrt{3}}+(\frac{1}{\sqrt{3}}-\frac{2}{3})Lt^{2k},\qquad\,\,\, x_{F_{k+1}}<x<x_{C_k}$\\ \hline
$\overline{D_kF_{k+1}}$  &$y=h_k^{DF}(x):=\sqrt{3}x+(1-\frac{2}{\sqrt{3}})Lt^{2k},\qquad\,\,\,\, x_{D_{k}}<x<x_{F_{k+1}}$\\\hline
$\overline{D_kE_{k+1}}$  &$x=h_k^{DE}(y):=(\frac{1}{2}-\frac{1}{\sqrt{3}})Lt^{2k},\qquad\qquad\quad\,\,\,\, y_{D_{k}}<y<y_{E_{k+1}}$\\\hline
$\overline{B_kE_{k+1}}$  &$y=h_k^{BE}(x):=-\sqrt{3}x+(1-\frac{2}{\sqrt{3}})Lt^{2k},\quad\,\,\,\,\,\, x_{B_{k}}<x<x_{E_{k+1}}$\\\hline
\end{tabular}
\caption{Sides of the polygons.}\label{1211091745}
\end{table}
We can now define the family of open sets whose union gives $\Omega$ (up to a set of zero Lebesgue measure):
\begin{itemize}
\item for $k=0$:
\begin{gather}
\omega_{A}:= \Bigl\{(x,y)\in\Omega: y<g_0^{BA}(x), x<g_0^{CA}(y), y>h_0^{BA}(x), y>h_0^{AC}(x) \Bigr\},\nonumber \\
\omega_{B}:= \Bigl\{(x,y)\in\Omega:y>g_0^{BA}(x), y>g_0^{EB}(x)\Bigr\},\nonumber \\
\omega_{C}:= \Bigl\{(x,y)\in\Omega: y>g_0^{CF}(x), x>g_0^{CA}(y)\Bigr\},\nonumber \\
\omega_{D}:= \Bigl\{(x,y)\in\Omega:y<g_0^{ED}(x), y<g_0^{DF}(x)\Bigr\},\nonumber \\
\omega_{E}:= \Bigl\{(x,y)\in\Omega: y>g_0^{ED}(x), y<g_0^{EB}(x), y<h_0^{BE}(x), x<h_0^{DE}(y)\Bigr\},\nonumber \\
\omega_{F}:= \Bigl\{(x,y)\in\Omega: y<g_0^{CF}(x), y>g_0^{DF}(x), y<h_0^{FC}(x), y<h_0^{DF}(x)\Bigr\},\nonumber
\end{gather}
\item for $k\geq 1$:
\begin{gather}
\omega_{A_k}:= \Bigl\{(x,y)\in\Omega: y<g_k^{BA}(x), x<g_k^{CA}(y), y>h_k^{BA}(x), y>h_k^{AC}(x) \Bigr\},\nonumber \\
\omega_{B_k}:= \Bigl\{(x,y)\in\Omega: y>g_k^{BA}(x), y>g_k^{EB}(x),y<h^{BA}_{k-1}(x),y>h_{k-1}^{BE}  \Bigr\},\nonumber \\
\omega_{C_k}:= \Bigl\{(x,y)\in\Omega: y>g_k^{CF}(x), x>g_k^{CA}(y), y<h_{k-1}^{AC}(x),y>h^{FC}_{k-1}(x)\Bigr\},\nonumber \\
\omega_{D_k}:= \Bigl\{(x,y)\in\Omega: y<g_k^{ED}(x), y<g_k^{DF}(x),x>h_{k-1}^{DE}(x),y>h_{k-1}^{DF}(x)\Bigr\},\nonumber \\
\omega_{E_k}:= \Bigl\{(x,y)\in\Omega: y>g_k^{ED}(x), y<g_k^{EB}(x), y<h_k^{BE}(x), x<h_k^{DE}(y)\Bigr\},\nonumber\\
\omega_{F_k}:= \Bigl\{(x,y)\in\Omega: y<g_k^{CF}(x), y>g_k^{DF}(x), y<h_k^{FC}(x), y<h_k^{DF}(x)\Bigr\},\nonumber
\end{gather}
\end{itemize}
With the exception of $\omega_B, \omega_C$ and $\omega_D$, the sets defined above are rhomboids (see Fig. \ref{1111121420} and \ref{1307081558}). All these rhomboids are similar to $\omega_A$, in the sense that they can be obtained from rescaling and rotation of $\omega_A$.
Indeed, for $k\geq 1$, the sets $\omega_{A_k}, \omega_{E_k}$ and $\omega_{F_k}$ are obtained by contracting each side of $\omega_{A}, \omega_{E}$ and $\omega_{F}$ respectively, by a factor $t^{2k}<1$. 
Moreover, by a $\frac{2\pi}{3}$ clockwise rotation of $\omega_A$ we obtain $\omega_F$ and by applying another $\frac{2\pi}{3}$ clockwise rotation we obtain $\omega_E$.
A similar relation holds for $\omega_{B_k},\omega_{C_k}$ and $\omega_{D_k}$ which are obtained by contracting each side of $\omega_{B_1},\omega_{C_1}$ and $\omega_{D_1}$ respectively by a factor of $t^{2k-2}$ (here with $k\geq 2$). Then, by a $\frac{2\pi}{3}$ clockwise rotation of $\omega_{B_1}$ we obtain $\omega_{C_1}$ and by applying another $\frac{2\pi}{3}$ clockwise rotation we obtain $\omega_{D_1}$.
Finally, by applying a counter-clockwise rotation of $\frac{\pi}{3}$ and then a contraction of each side by $t$ to $\omega_A$ we obtain $\omega_{B_1}$.
The measure of the four internal angles of $\omega_{A}$ are respectively $\frac{5\pi}{6},\frac{\pi}{2},\frac{\pi}{2},\frac{\pi}{6}$. By similarity, the internal angles of all the rhomboids defined above coincide with those of $\omega_A$. From the equations of the segments reported in Table \ref{1211091745} it is also possible to compute the normals to each side of the rhomboid. 
By denoting with $\calL^2$ the (two-dimensional) Lebesgue measure, it follows that 
\begin{eqnarray}\label{1307111200}
\calL^2(\omega_{A_k})=\calL^2(\omega_{E_k})=\calL^2(\omega_{F_k})
\textrm{  and } \calL^2(\omega_{B_k})=\calL^2(\omega_{C_k})=\calL^2(\omega_{D_k}) \quad \forall k\geq 0.
\end{eqnarray}
Since for $k\geq 1$ $\omega_{B_k}$ is obtained (up to a rotation) by contracting $\omega_{A_{k-1}}$ by a factor $t<1$, we have $\calL^2(\omega_{B_k})=t^{2}\calL^2(\omega_{A_{k-1}})$ and, since $\omega_{A_{k}}$ is obtained by contracting each side of $\omega_{A_k-1}$ by a factor $t^2$, we have $\calL^2(\omega_{A_{k}})=t^{4}\calL^2(\omega_{A_k-1})$.
Summarizing:
\begin{eqnarray}\label{1307081703}
\calL^2(\omega_{A_k})=\calL^2(\omega_{E_k})=\calL^2(\omega_{F_k})=L^2 t^{4k}t,\, \calL^2(\omega_{B_k})=\calL^2(\omega_{C_k})=\calL^2(\omega_{D_k})=L^2 t^{4k}t^{-1}.
\end{eqnarray}
If we define
$$
\omega_k:=\omega_{A_k}\cup\omega_{B_k}\cup\omega_{C_k}\cup\omega_{D_k}\cup\omega_{E_k}\cup\omega_{F_k}
$$
we have $\bigcup\omega_k\subset\Omega$. Importantly, as the remainder in this inclusion is represented by elements of $\calL^2$-measure equal to zero we have
\begin{eqnarray}\label{1306042258}
\calL^2(\Omega)=\sum_{k=0}^{\infty}\calL^2(\omega_k).
\end{eqnarray}
Notice that, for any $ k\geq 0$, the origin does not belong to any of the sets $\omega_k$. Furthermore, each neighborhood of the origin contains an infinite number of sets $\omega_k$. This fact has deep consequence on the property of the microstructure and it is at the origin of the nesting of variants close to the origin.

\subsection{Construction of the microstructure}\label{1211121543}

As typical when dealing with microstructure, the definition of the map $\uu$ is based on the introduction of a piecewise-constant tensor field $\HH(x,y)$ such that $\nabla\uu(x,y)=\HH(x,y)$. 
It is widely known that not every matrix field $\HH(x,y)$ is the gradient of a continuous vector field. In fact, this condition reveals as a geometrical compatibility constraint on $\HH(x,y)$ which is known as Hadamard jump condition. Simply, consider a smooth line $S$ having normal $\nu$. If $\uu$ is continuous at both sides of $S$ with limits $\nabla\uu^+(x,y)=\HH^+$ and $\nabla\uu^-(x,y)=\HH^-$ at each point $(x,y)\in S$, from above and below $S$ respectively, then equating the tangential derivatives leads to the Hadamard jump condition
\begin{eqnarray}\label{1306142056}
\HH^+-\HH^-=\aa\otimes \nn
\end{eqnarray}
with $\aa,\nn\in\RRR^2$ and $|\nn|=1$. The matrices $\HH^+$ and $\HH^-$ in (\ref{1306142056}) are said to be rank-one connected. Notice that for $2\times 2 $ matrices this is equivalent to the condition $\det(\HH^+-\HH^-)=0$.
The situation in which we have a map $\uu$ with $\nabla\uu$ piecewise constant and equal to either $\HH^+$ or $\HH^-$ in alternating bands bounded by lines perpendicular to $\nn$ is relevant in materials science. Such maps are called simple laminates \cite{Ball87}, \cite{Muller}.
With no further hypotheses on the domain, condition $(\ref{1306142056})$ is a necessary condition for the existence of a simple laminate $\uu$ such that $\nabla\uu=\HH^+,\HH^-$. 
In the case in which we have more interfaces, condition $(\ref{1306142056})$ has to be verified at each interface across which we have two  crystal phases corresponding to different deformation gradients.

Besides the Hadamard jump condition, which can be regarded as a \textit{geometric} compatibility condition, the construction of the map $\HH(x,y)$ has to satisfy an \textit{energetic} condition. In fact, this microstructure is zero-energy in the following sense.
Recall that $\varepsilon_2$ and $\varepsilon_3$ are the components of the deviatoric part of $\EE$ and $\varepsilon_1$ is the unique component of the spherical part of $\EE$ introduced in (\ref{1307041431}). If we assume that $\psi_{L}$ has three distinct minima corresponding to the three strains (\ref{1305131529}) then
\begin{eqnarray}\label{1307070104}
\int_{\Omega}\psi_{L}\bigl(\varepsilon_1,\varepsilon_2,\varepsilon_3\bigr) dxdy=0
\end{eqnarray}
if and only if
\begin{eqnarray}\label{star}
\bigl(\nabla\uu(x,y)\bigr)_{sym}\in\{\EE_1,\EE_2,\EE_3\} \,\,\textrm{for a.e. } (x,y) \in \Omega.
\end{eqnarray}
Therefore, the differential inclusion (\ref{star}) guarantees that the map $\uu$ is a minimizer for $\int_{\Omega}\psi_{L}(\varepsilon_2,\varepsilon_3)dxdy$.
Since $\HH_{sym}$ has to be equal to either $\EE_1$, $\EE_2$ or $\EE_3$ for almost every $(x,y)\in\Omega$, what is left to determine is the skew part of $\HH(x,y)$.

For the reader's convenience, we illustrate with an example how the map $\HH$ can fulfill the geometric constraint and the energetic constraint simultaneously.
Here we make use of some notation that will be further defined and used extensively in the next section.
Consider the situation of Fig. \ref{1111121420}-LEFT, and more precisely the matching between the domains $\omega_A$ and $\omega_B$ . We construct a piecewise constant tensor field $\HH$ defined as 
\begin{eqnarray}\label{1306161632}
\HH=\left\{
\begin{array}{cccc}
\HH_A &\textrm{on} &\omega_A\\
\HH_B &\textrm{on} &\omega_B
\end{array} \right.
\end{eqnarray}
where the matrices $\HH_A$ and $\HH_B$ are to be determined
in a way such that $(\HH_{A})_{sym}=\EE_3$ and $(\HH_{B})_{sym}=\EE_1$ and $(\HH_{A}-\HH_{B})=\alpha\otimes\nu$, for some $\alpha,\nu\in\RRR^2$ with $|\nu|=1$.
We recall that $\EE_1$ and $\EE_3$ are defined in (\ref{1305131529}).
Let us write $\HH_A=\EE_3$ and ${\HH_B}=\EE_1+\WW$. Here $\WW$ is an unknown skew-symmetric matrix written in the form,
\begin{eqnarray}\label{}
\WW= \epsilon\left(
\begin{array}{ccc}
0 & -w  \\
w & 0
\end{array} \right),
\quad w\in\RRR.
\end{eqnarray}
The condition $(\ref{1306142056})$ applied to $\HH_A$ and ${\HH_B}$ reads
$(\EE_3-\EE_1-\WW)=\aa\otimes \nn$
or, componentwise,
\begin{eqnarray}\label{1306161632}
\left\{
\begin{array}{cccc}
-\frac{3}{2}\epsilon &=a_1 n_1\\
(-\frac{\sqrt{3}}{2}+w)\epsilon &=a_1 n_2\\
(-\frac{\sqrt{3}}{2}-w)\epsilon &=a_2 n_1\\
\frac{3}{2}\epsilon &=a_2 n_2\\
\end{array} \right.
\end{eqnarray}
with $n_1^2+n_2^2=1$. The system  of equations (\ref{1306161632}) has several solutions. Let us choose the solution with $n_1^2=\frac{3}{4}$, $w=\sqrt{3}$ and vectors $\alpha:=(a_1,a_2)=\epsilon(\sqrt{3},3)$, $\nu:=(n_1,n_2)=(-\frac{\sqrt{3}}{2},\frac{1}{2})$. The vector $\nu$ is perpendicular to $\overline{BA}$, the interface separating $\omega_A$ and $\omega_B$. Thus, we obtain,
\begin{eqnarray}\label{1306201635}
\HH_{B}=\EE_1+
\epsilon\left(
\begin{array}{ccc}
0 & -\sqrt{3}  \\
\sqrt{3} & 0
\end{array} \right).
\end{eqnarray}
 The fact that the vector $\nu$ is exactly orthogonal to the interface separating $\omega_{A}$ and $\omega_{B}$ guarantees the existence of a continuous map $\uu$ such that its gradient coincides with $\HH$ on $\omega_{A}\cup\omega_{B}$.\footnote{We note that the geometric constraint has been implicitly taken into account in the construction of the star-shaped geometry of the previous section.}

We now focus on the matching between the domains $\omega_{A}$ and $\omega_{B_1}$. As the strains associated with the domains $\omega_{B}$ and   $\omega_{B_1}$  are equal, the matching equations to be solved are the same than for the matching between $\omega_{A}$ and $\omega_{B}$, that is, Eqs. (\ref{1306161632}). 
In this case, however, we consider a different solution,  namely, 
$n_1^2=\frac{1}{4}$, $w=-\sqrt{3}$ and vectors $\beta:=(a_1,a_2)=\epsilon(-3,\sqrt{3})$, $\mu:=(n_1,n_2)=(\frac{1}{2},\frac{\sqrt{3}}{2})$, where the vector $\mu$ is perpendicular to $\overline{B_1A}$, the interface separating $\omega_A$ and $\omega_{B_1}$. Thus, in this case we obtain, 
\begin{eqnarray}\label{1306201637}
\HH_{B_1}=\EE_1+
\epsilon\left(
\begin{array}{ccc}
0 & \sqrt{3}  \\
-\sqrt{3} & 0
\end{array} \right).
\end{eqnarray}
Notice that the symmetric part of $\HH_{B_1}$ (as well as the symmetric part of $\HH_B$ as defined in Eq. (\ref{1306201635})) is precisely $\EE_1$. 
The fact that $\mu$ is orthogonal to the interface separating $\omega_A$ and $\omega_{B_1}$ is enough to guarantee that there exists a function $\widehat{\uu}$ such that $\nabla\widehat{\uu}=\HH_A$ on $\omega_A$ and $\nabla\widehat{\uu}=\HH_{B_1}$ on $\omega_{B_1}$. Of course, if $\uu\equiv\widehat{\uu}$ over $\omega_A$ one can define a continuous function defined over the larger set $\omega_A\cup\omega_B\cup\omega_{B_1}$.

\noindent Even though the example above is restricted to three subsets of $\Omega$, it suggests that this technique requires a careful matching of matrices and depends heavily on the geometry of the problem. Whether it is possible to further iterate this construction over a larger domain is not a trivial question. 
Here one of the mathematical difficulties lies in the fact that this microstructure becomes finer and finer close to the center of the tripole-star pattern. 
%
%
Indeed, as $k\to\infty$, the size of the sets $\omega_k$ is reduced while the number of interfaces increases. Since one must be able to prove that the Hadamard jump condition is satisfied at each interface, it is not a priori clear if the construction of a continuous map $\uu$ can be extended asymptotically for each hierarchy of laminates. In other words, the point $(0,0)$ may act as a topological singularity. Below we show that this construction can actually be extended to the whole $\Omega$.
Additionally, in accordance with experimental observations, our construction involves a mixture of variants occurring with the same phase fraction. By defining
\begin{eqnarray}\label{1307111218}
\omega_{\EE_i}:= \Bigl\{(x,y)\in\Omega:\bigl(\nabla\uu(x,y)\bigr)_{sym}=\EE_i\Bigr\},\quad i=1,2,3,
\end{eqnarray}
the sets of points $(x,y)$ in $\Omega$ such that $(\nabla\uu)_{sym}=\EE_i$, we have that $\omega_{\EE_1}$, $\omega_{\EE_2}$, and $\omega_{\EE_3}$ have the same area.


In the remaining part of this section we construct the function $\uu$ from scratch instead of focusing on the abstract properties for which such a construction can be achieved. This approach has, of course, an intrinsic limitation in the sense that it is not clear if it can be extended to treat more general cases of non-periodic microstructures.
The construction of $\uu$ will serve as the proof of the following theorem, which is the central result of this paper.

\begin{Theorem}\label{THM}
Let $L$ and $\epsilon$ be any two positive real numbers.
Let $\Omega=B_{R}(0,0)$ with $R=L\sqrt{\frac{2}{3}+\frac{1}{\sqrt{3}}}$, let $\EE_i$ be defined as in $(\ref{1305131529})$
and $\omega_{\EE_i}$ as defined in {\rm  (\ref{1307111218})}. Then there exists a map $\uu:\Omega\to\RRR^2$ such that 
\begin{enumerate}
\item[1)] $\uu$ is continuous in $\overline{\Omega}$,
\item[2)] $\bigl(\nabla\uu(x,y)\bigr)_{sym}\in \{\EE_1,\EE_2,\EE_3\} \,\,\, \textrm{for a.e. (x,y)}\,\, \textrm{in} \,\,\Omega,$
\item[3)] $\calL^2\bigl(\omega_{\EE_i}\bigr)=\frac{\pi}{3} R^2, \quad i=1,2,3,$
\item[4)] $\uu$ lies in (an equivalence class in) $W^{1,p}(\Omega,\RRR^2), \,\, 1\leq p<\infty$. More precisely, $\uu$ is not Lipschitz continuous.

\end{enumerate}
 
\end{Theorem}
 
\begin{Remark}
Proof of Theorem \ref{THM} lies in a tedious but straightforward series of algebraic computations. To keep the presentation engaging, we give here a sketch of the argument required.
%
%
To construct the function $\uu$, we must ensure that the geometrical compatibility condition holds at each interface.
Indeed, consider the two rank-one connected matrices $\HH_A$ and $\HH_{B}$, introduced in the example above. Since compatibility holds at the interface $\overline{BA}$, upon integration of $\HH_{A}$ and $\HH_{B}$ we can define, up to constants of integration, the two functions
\begin{gather}
\epsilon\vv_{A}(x,y) =\epsilon\left\{
\begin{array}{cccc}
-\frac{x}{2}-\frac{\sqrt{3}}{2}y+c_{A,1}\\
-x\frac{\sqrt{3}}{2}+\frac{y}{2}+c_{A,2}\\
\end{array} \right.\nonumber\qquad
\epsilon\vv_{B}(x,y) =\epsilon\left\{
\begin{array}{cccc}
x-\sqrt{3} y+c_{B,1}\\
\sqrt{3}x-y+c_{B,2} \\
\end{array} \right.\nonumber
\end{gather}
over $\omega_A$ and $\omega_B$ respectively. By requiring continuity at the interface
we determine $c_{B,1}$ and $c_{B,2}$, thus leaving $c_{A_1}$ and $c_{A,2}$ unknown. 
Similarly, matrices
$\HH_C$, $\HH_D$, $\HH_E$ and $\HH_F$
can be defined by requiring their symmetric part be equal to either $\EE_1, \EE_2$ or $\EE_3$ 
and by requiring them to be pairwise rank-one connected thus
fixing their skew parts.
%
%
%
Upon integration of these matrices one defines vectors $\vv_C$, $\vv_D$, $\vv_E$ and $\vv_F$ over $\omega_C$, $\omega_D$, $\omega_E$ and $\omega_F$ respectively. Notice that the constants of integration will again be determined by matching the vectors at the interfaces (except $c_{A,1}$ and $c_{A,2}$ which are undetermined). We have then constructed a piece-wise affine vector field which is continuous over $\omega_0$ and which we denote with $\vv_0$.
An analogous procedure over $\omega_k$ yields a piecewise affine vector field which we denote with $\vv_k$.
%
As for $\vv_0$, each function $\vv_k$ with $k\geq 1$ is defined up to constants of integration.
Finally, the introduction of the globally continuous map $\uu$ follows after the computation of the remaining constants of integration (up to two which remain, correctly, unknown) by matching the maps $\vv_{k}$ and $\vv_{k+1}$ for every $k\geq 0$. This is in fact possible because the geometric compatibility condition holds also for each interfaces between $\omega_k$ and $\omega_{k+1}$.

\end{Remark}

\paragraph{Proof of Theorem \ref{THM}.}
In what follows we make an extensive use of the notation introduced in Paragraph \ref{1003151615}. Let us define
\begin{eqnarray}\label{1307072020}
\widetilde{\uu}(x,y)=
\epsilon\vv(x,y)+\epsilon\vv_o(x,y),
\end{eqnarray}
where $\vv:\bigcup_{k=0}^{\infty}\omega_k\to\RRR^2$ is defined as $\vv=\vv_k$ on $\omega_k$ and $\vv_k:\omega_k\to\RRR^2$ is defined as
\begin{eqnarray}
\vv_k(x,y) :=\left\{
\begin{array}{cccc}
\vv_{A_k} &\textrm{on}\quad\omega_{A_k}\\
\vv_{B_k} &\textrm{on}\quad\omega_{B_k}\\
\vv_{C_k} &\textrm{on}\quad\omega_{C_k}\\
\vv_{D_k} &\textrm{on}\quad\omega_{D_k}\\
\vv_{E_k} &\textrm{on}\quad\omega_{E_k}\\
\vv_{F_k} &\textrm{on}\quad\omega_{F_k}\\
\end{array} \right.
\quad k\geq 0
\end{eqnarray}
with
\begin{eqnarray}
\vv_{A_k}(x,y)&:=&\displaystyle{\begin{dcases}
-\frac{x}{2}+\Bigl(2k-\frac{1}{2}\Bigr)\sqrt{3}y \\
-\sqrt{3}x\Bigl(\frac{1}{2}+2k\Bigr)+\frac{y}{2}+Lt^{2k}
\end{dcases}} \nonumber\\
\vv_{B_k}(x,y)&:=&\displaystyle{\begin{dcases}
x+(-\sqrt{3}+2k\sqrt{3} )y+L\frac{t^{2k}}{2}   \\
x(\sqrt{3}-2k\sqrt{3})-y+\frac{t^{2k}}{2}(2+\sqrt{3})L
\end{dcases}} \nonumber\\
\vv_{C_k}(x,y)&:=&\displaystyle{\begin{dcases}
 -\frac{x}{2}+\sqrt{3}y\Bigl(2k-\frac{1}{2}\Bigr) \\
\frac{3\sqrt{3}}{2}x-2k\sqrt{3}x+\frac{y}{2} 
\end{dcases}} \nonumber\\
\vv_{D_k}(x,y)&:=&\displaystyle{\begin{dcases}
-\frac{x}{2}+2k\sqrt{3}y-\frac{3\sqrt{3}}{2}y-t^{2k}\Bigl(\frac{1}{2}+\frac{\sqrt{3}}{2}\Bigr)L \\
\Bigl(\frac{\sqrt{3}}{2}-2k\sqrt{3}\Bigr)x+\frac{y}{2}+t^{2k}\Bigl(\frac{1+\sqrt{3}}{2}\Bigr)L 
\end{dcases}} \nonumber
\end{eqnarray}
\begin{eqnarray}
\vv_{E_k}(x,y)&:=&\displaystyle{\begin{dcases}
-\frac{x}{2}+\Bigl(\frac{1}{2}+2k\Bigr)\sqrt{3}y+t^{2k}\Bigl(\frac{1}{2}-\frac{\sqrt{3}}{2}\Bigr)L \\
\Bigl(\frac{1}{2}-2k\Bigr)\sqrt{3}x+\frac{y}{2}+t^{2k}\Bigl(\frac{1}{2}+\frac{\sqrt{3}}{2}\Bigr)L 
\end{dcases}} \nonumber\\
\vv_{F_k}(x,y)&:=&\displaystyle{\begin{dcases}
x+2\sqrt{3}ky-\frac{1}{2}t^{2k}L\\
-2k\sqrt{3}x-y+t^{2k}\frac{\sqrt{3}}{2}L 
\end{dcases}} \nonumber
\end{eqnarray}
and $\vv_o:\bigcup_{k=0}^{\infty}\omega_k\to\RRR^2$ is defined as $\vv_o=\vv_{ok}$ on $\omega_k$ and $\vv_{ok}:\omega_k\to\RRR^2$ is defined as
\begin{eqnarray}
\textrm{for } k=0:\qquad \vv_{o0}(x,y)&:=&\displaystyle{\begin{dcases}
0\\
0
\end{dcases}} \nonumber\\
\textrm{for } k\geq 1: \qquad  \vv_{ok}(x,y)&:=&\displaystyle{\begin{dcases}
L(\sqrt{3}-2)\sum_{j=1}^k t^{2j-2}  \\
L\sum_{j=1}^k t^{2j-2}.
\end{dcases}} \nonumber
\end{eqnarray}
By construction, the function $\widetilde{\uu}$ is defined at almost every point in $\Omega$.
First of all, we show that there exists a continuous extension ${\uu}$ of $\widetilde{\uu}$ defined at each point $(x,y)$ of $\Omega$. Notice that, on each subset $\omega_{A_k}$,$\omega_{B_k}$, $\omega_{C_k}$, $\omega_{D_k}$, $\omega_{E_k}$ and $\omega_{F_k}$ the map $\widetilde{\uu}$ is affine and therefore continuous. We now check that $\widetilde{\uu}$ has indeed continuous limit across each interface. In what follows we compute the limit $\widetilde{\uu}^+, \widetilde{\uu}^-$ at each interface. 

To begin, consider all the interfaces inside $\omega_k$. Let $k\geq 0$:
\begin{enumerate}
\item Interface $\overline{B_k A_k}$  
\begin{eqnarray}
\widetilde{\uu}^-(x,y):=\lim_{(z_x,z_y)\to (x,y)^-}\widetilde{\uu} (z_x,z_y)=\lim_{(z_x,z_y)\to (x,y)}\epsilon(\vv_{B_k}+\vv_{ok})(z_x,z_y)=\qquad\qquad\qquad\nonumber\\
\qquad\qquad=\lim_{(z_x,z_y)\to (x,y)}\epsilon(\vv_{A_k}+\vv_{ok})(z_x,z_y)=\lim_{(z_x,z_y)\to (x,y)^+}\widetilde{\uu} (x,y)=:\widetilde{\uu}^+(x,y)=\nonumber\\
=\epsilon\vv_{ok}+\epsilon\left\{
\begin{array}{cccc}
 -2x+6kx+L(2k-\frac{1}{2})t^{2k} \\
-2\sqrt{3}kx +\frac{L}{2\sqrt{3}}t^{2k}+Lt^{2k} 
\end{array} \right. \nonumber
\end{eqnarray}
\item Interface $\overline{E_k B_k}$
\begin{eqnarray}
\lim_{(z_x,z_y)\to (x,y)}\epsilon(\vv_{B_k}+\vv_{ok})(z_x,z_y)=\lim_{(z_x,z_y)\to (x,y)}\epsilon(\vv_{E_k}+\vv_{ok})(z_x,z_y)=\qquad\qquad\qquad\nonumber\\
\qquad\qquad\qquad=
\epsilon\vv_{ok}+\epsilon\left\{
\begin{array}{cccc}
2kx+Lt^{2k}\frac{2k}{\sqrt{3}}+Lt^{2k}\frac{1}{2\sqrt{3}}(\sqrt{3}-2)  \\
\frac{2}{\sqrt{3}}x-2k\sqrt{3}x+\frac{L}{2}\sqrt{3}t^{2k}+\frac{2}{3}t^{2k}L 
\nonumber
\end{array} \right.
\end{eqnarray}
\item Interface $\overline{E_k D_k}$
\begin{eqnarray}
\lim_{(z_x,z_y)\to (x,y)}\epsilon(\vv_{D_k}+\vv_{ok})(z_x,z_y)=\lim_{(z_x,z_y)\to (x,y)}\epsilon(\vv_{E_k}+\vv_{ok})(z_x,z_y)=\qquad\qquad\quad\nonumber\\
\qquad\qquad\qquad=\epsilon\vv_{ok}+\epsilon\left\{
\begin{array}{cccc}
-\frac{x}{2}-Lt^{2k}k -t^{2k}(\frac{\sqrt{3}}{2}-\frac{1}{4})L \\
\frac{\sqrt{3}}{2}x-2k\sqrt{3}x+(-\frac{L}{4\sqrt{3}}t^{2k})+t^{2k}(\frac{1+\sqrt{3}}{2})L \\
\end{array} \right.\nonumber
\end{eqnarray}
\item Interface $\overline{D_k F_k}$
\begin{eqnarray}
\lim_{(z_x,z_y)\to (x,y)}\epsilon(\vv_{D_k}+\vv_{ok})(z_x,z_y)=\lim_{(z_x,z_y)\to (x,y)}\epsilon(\vv_{F_k}+\vv_{ok})(z_x,z_y)=\quad\qquad\qquad\nonumber\\
\qquad\qquad\qquad=\epsilon\vv_{ok}+\epsilon\left\{
\begin{array}{cccc}
x-2kx-\frac{L}{\sqrt{3}}t^{2k}2k -t^{2k}\frac{1}{2}L \\
\frac{1}{\sqrt{3}}x-2k\sqrt{3}x+\frac{t^{2k}}{3}L+L\frac{\sqrt{3}t^{2k}}{2} \\
\end{array} \right.\nonumber
\end{eqnarray}
\item Interface $\overline{C_k F_k}$
\begin{eqnarray}
\lim_{(z_x,z_y)\to (x,y)}\epsilon(\vv_{C_k}+\vv_{ok})(z_x,z_y)=\lim_{(z_x,z_y)\to (x,y)}\epsilon(\vv_{F_k}+\vv_{ok})(z_x,z_y)=\quad\qquad\qquad\nonumber\\
\qquad\qquad\qquad=\epsilon\vv_{ok}+\epsilon\left\{
\begin{array}{cccc}
 x-6kx-\frac{Lt^{2k}}{2}+2kLt^{2k} \\
-2k\sqrt{3}x+\sqrt{3}x+\frac{L}{2\sqrt{3}}t^{2k} \\
\end{array} \right.\nonumber
\end{eqnarray}
\item Interface $\overline{C_k A_k}$
\begin{eqnarray}
\lim_{(z_x,z_y)\to (x,y)}\epsilon(\vv_{C_k}+\vv_{ok})(z_x,z_y)=\lim_{(z_x,z_y)\to (x,y)}\epsilon(\vv_{A_k}+\vv_{ok})(z_x,z_y)=\qquad\qquad\qquad\nonumber\\
\qquad\qquad\qquad=\epsilon\vv_{ok}+\epsilon\left\{
\begin{array}{cccc}
-\frac{L}{4\sqrt{3}}t^{2k}+2k\sqrt{3}y-\frac{\sqrt{3}}{2}y \\
-kLt^{2k}+\frac{y}{2}+\frac{3}{4}Lt^{2k} .\\
\end{array} \right.\nonumber
\end{eqnarray}
\end{enumerate}
Let us now consider the interfaces between $\omega_k$ and $\omega_{k+1}$, for any $k\geq 0$:
\begin{enumerate}
\item Interface $\overline{B_k A_{k+1}}$
\begin{eqnarray}
\widetilde{\uu}^-(x,y):=\lim_{(z_x,z_y)\to (x,y)^-}\widetilde{\uu} (z_x,z_y)=\lim_{(z_x,z_y)\to (x,y)}\epsilon(\vv_{B_{k+1}}+\vv_{ok+1})(z_x,z_y)=\qquad\quad\nonumber\\
\qquad\quad=\lim_{(z_x,z_y)\to (x,y)}\epsilon(\vv_{A_{k}}+\vv_{ok})(z_x,z_y)=\lim_{(z_x,z_y)\to (x,y)^+}\widetilde{\uu}(x,y)=:\widetilde{\uu}^+(x,y)=\nonumber\\
=\epsilon\vv_{ok}+\epsilon\left\{
\begin{array}{cccc}
-2kx+2kLt^{2k}(\frac{2}{\sqrt{3}}-1)+t^{2k}L(\frac{1}{2}-\frac{1}{\sqrt{3}})\\
-\frac{2x}{\sqrt{3}}-2\sqrt{3}kx+Lt^{2k}(\frac{4}{3}-\frac{1}{2\sqrt{3}})
\\
\end{array} \right.\nonumber
\end{eqnarray}
\item Interface $\overline{B_k E_{k+1} }$
\begin{eqnarray}
\lim_{(z_x,z_y)\to (x,y)}\epsilon(\vv_{B_{k+1}}+\vv_{ok+1})(z_x,z_y)=\lim_{(z_x,z_y)\to (x,y)}\epsilon(\vv_{E_k}+\vv_{ok})(z_x,z_y)=\nonumber\\
=\epsilon\vv_{ok}+\epsilon\left\{
\begin{array}{cccc}
-2x-6kx+Lt^{2k}(-\frac{1}{2}-4k+2k\sqrt{3})
\\
-2kx\sqrt{3}+Lt^{2k}(1+\frac{1}{2\sqrt{3}})
\\
\end{array} \right.\nonumber
\end{eqnarray}
\item Interface $\overline{D_k E_{k+1} }$
\begin{eqnarray}
\lim_{(z_x,z_y)\to (x,y)}\epsilon(\vv_{D_{k+1}}+\vv_{ok+1})(z_x,z_y)=\lim_{(z_x,z_y)\to (x,y)}\epsilon(\vv_{E_k}+\vv_{ok})(z_x,z_y)=\qquad\nonumber\\=\epsilon\vv_{ok}+\epsilon\left\{
\begin{array}{cccc}
\sqrt{3}(\frac{1}{2}+2k)y+Lt^{2k}(\frac{1}{4}-\frac{1}{\sqrt{3}})\\
\frac{y}{2}+Lt^{2k}(2k-k\sqrt{3}+\frac{3}{4}\sqrt{3})\\
\end{array} \right.\nonumber
\end{eqnarray}
\item Interface  $\overline{D_k F_{k+1}}$
\begin{eqnarray}
\lim_{(z_x,z_y)\to (x,y)}\epsilon(\vv_{D_{k+1}}+\vv_{ok+1})(z_x,z_y)=\lim_{(z_x,z_y)\to (x,y)}\epsilon(\vv_{F_k}+\vv_{ok})(z_x,z_y)=\qquad\nonumber\\
\qquad=\epsilon\vv_{ok}+\epsilon\left\{
\begin{array}{cccc}
 x+6kx-\frac{1}{2}Lt^{2k} +2k\sqrt{3}Lt^{2k}-4kLt^{2k}\\
x(-\sqrt{3}-2k\sqrt{3})+Lt^{2k}(-1+\frac{7}{2\sqrt{3}})\\
\end{array} \right.\nonumber
\end{eqnarray}
\item Interface $\overline{F_{k+1}C_k}$
\begin{eqnarray}
\lim_{(z_x,z_y)\to (x,y)}\epsilon(\vv_{C_{k+1}}+\vv_{ok+1})(z_x,z_y)=\lim_{(z_x,z_y)\to (x,y)}\epsilon(\vv_{F_k}+\vv_{ok})(z_x,z_y)=\qquad\nonumber\\
=\epsilon\vv_{ok}+\epsilon\left\{
\begin{array}{cccc}
x(1+2k)+Lt^{2k}(2k-\frac{4k}{\sqrt{3}}-\frac{1}{2}) \\
x(-2k\sqrt{3}-\frac{1}{\sqrt{3}})+Lt^{2k}(\frac{2}{3}+\frac{1}{2\sqrt{3}}) \\
\end{array} \right.\nonumber
\end{eqnarray}
\item Interface $\overline{A_{k+1} C_k}$
\begin{eqnarray}
\lim_{(z_x,z_y)\to (x,y)}\epsilon(\vv_{C_{k+1}}+\vv_{ok+1})(z_x,z_y)=\lim_{(z_x,z_y)\to (x,y)}\epsilon(\vv_{A_k}+\vv_{ok})(z_x,z_y)=\qquad\nonumber\\=\epsilon\vv_{ok}+\epsilon\left\{
\begin{array}{cccc}
-\frac{x}{2}+Lt^{2k}(\frac{\sqrt{3}}{4}-\frac{1}{2}+2k-\sqrt{3}k) \\
-\frac{\sqrt{3}}{2}x-2k\sqrt{3}x+Lt^{2k}(\frac{1}{2\sqrt{3}}+\frac{3}{4}) .\\
\end{array} \right.\nonumber
\end{eqnarray}
\end{enumerate}
The above calculations show it is possible to extend $\widetilde{\uu}$ to a function still denoted by $\widetilde{\uu}$ by continuity across each interface. Notice that the function $\widetilde{\uu}$ is continuous also at the vertices $A_k$, $B_k$, $C_k$, $D_k$, $E_k$ and $F_k$.
The new map $\widetilde{\uu}$ is now defined at each point of $\Omega/(0,0)$ and, for each $\delta>0$, $\widetilde{\uu}$ is continuous on the set $\Omega/\overline{B_{\delta}(0,0)}$. Therefore, in order to extend $\widetilde{\uu}$ to a continuous function on the whole $\Omega$, we have to clarify the behavior of $\widetilde{\uu}$ close to the point $(0,0)$.
Precisely, we prove that
\begin{eqnarray}\label{1306252145}
\lim_{(z_x,z_y)\to(0,0)}\widetilde{\uu}(z_x,z_y)=\epsilon L(\sqrt{3}-2,1)\frac{1}{1-t^2}.
\end{eqnarray} 
%
%
Notice that it is enough to show that the supremum and the infimum of the components of $\widetilde{\uu}$ on $\omega_k$ generate sequences $\{\MM_k\}=\{M_{1,k},M_{2,k}\}$ and $\{\mm_k\}=\{m_{1,k},m_{2,k}\}$ both converging to $\epsilon L(\sqrt{3}-2,1)\frac{1}{1-t^2}$ as $k\to\infty$. This, combined with a standard contradiction argument guarantees that  (\ref{1306252145}) holds for every sequence of vectors converging to $(0,0)$. 
%
%
Let us denote with $\widetilde{u_i}(x,y)$, with $i=1,2$, the components of $\widetilde{\uu}(x,y)$. We compute the supremum of $\widetilde{u_i}$ by splitting the computation over the subdomains composing $\omega_k$ yielding
\begin{eqnarray}\label{13060121244}
M_{i,k}:=\sup_{(x,y)\in\omega_{k}}\widetilde{u_i}(x,y)=\max\Bigl\{\Bigl(\sup_{(x,y)\in\omega_{A_k}}\widetilde{u_i}(x,y)\Bigr),\Bigl(\sup_{(x,y)\in\omega_{B_k}}\widetilde{u_i}(x,y)\Bigr),\nonumber\\\Bigl(\sup_{(x,y)\in\omega_{C_k}}\widetilde{u_i}(x,y)\Bigr), \Bigl(\sup_{(x,y)\in\omega_{D_k}}\widetilde{u_i}(x,y)\Bigr),
\Bigl(\sup_{(x,y)\in\omega_{E_k}}\widetilde{u_i}(x,y)\Bigr),\Bigl(\sup_{(x,y)\in\omega_{F_k}}\widetilde{u_i}(x,y)\Bigr)\Bigr\}.
\end{eqnarray}
Since $\widetilde{\uu}$ is affine over each of the subdomains in (\ref{13060121244}), supremum and infimum are attained on the boundary (and, in particular, on the vertices). Therefore, focusing on the case of $\omega_{A_k}$ we can write
\begin{eqnarray}\label{13060123222}
\sup_{(x,y)\in\omega_{A_k}}\widetilde{u_i}(x,y)=\max\bigl\{\widetilde{u_i}(x_{A_k},y_{A_k}),\,\widetilde{u_i}(x_{B_k},y_{B_k}),\,\widetilde{u_i}(x_{C_k},y_{C_k}),\,\widetilde{u_i}(x_{A_{k+1}},y_{A_{k+1}})\bigr\}.
\end{eqnarray}
A straightforward computation yields
\begin{eqnarray}\label{1306252230}
\widetilde{\uu}(x_{A_k},y_{A_k})=\epsilon\left\{
\begin{array}{cccc}
Lt^{2k}\bigl[-\frac{1}{\sqrt{3}}-\frac{1}{2}+k(2+\sqrt{3})\bigr]+\bigl(\sqrt{3}-2\bigr)L\sum_{j=1}^k t^{2j-2}\\
Lt^{2k}\bigl[\frac{1}{2\sqrt{3}}+1-k \bigr]+L\sum_{j=1}^k t^{2j-2}
\end{array} \right.
\end{eqnarray}
\begin{eqnarray}
\widetilde{\uu}(x_{B_k},y_{B_k})=\epsilon\left\{
\begin{array}{cccc}
Lt^{2k}\bigl[ k(\sqrt{3}-1)+\frac{1}{2}-\frac{1}{\sqrt{3}}\bigr]+\bigl(\sqrt{3}-2\bigr)L\sum_{j=1}^k t^{2j-2}\\
Lt^{2k}\bigl[ k(\sqrt{3}-1)+1+\frac{1}{2\sqrt{3}}\bigr]+L\sum_{j=1}^k t^{2j-2}
\end{array} \right.
\end{eqnarray}
\begin{eqnarray}
\widetilde{\uu}(x_{C_k},y_{C_k})=\epsilon\left\{
\begin{array}{cccc}
Lt^{2k}\bigl[ k(2-\sqrt{3})+\frac{\sqrt{3}}{6}-\frac{1}{2}\bigr]+\bigl(\sqrt{3}-2\bigr)L\sum_{j=1}^k t^{2j-2}\\
Lt^{2k}\bigl[ -k +\frac{1}{2}+\frac{1}{2\sqrt{3}}\bigr]+L\sum_{j=1}^k t^{2j-2}\\
\end{array} \right.
\end{eqnarray}
\begin{eqnarray}\label{1306252234}
\widetilde{\uu}(x_{A_{k+1}},y_{A_{k+1}})=\epsilon\left\{
\begin{array}{cccc}
Lt^{2k+2}\bigl[-\frac{1}{\sqrt{3}}-\frac{1}{2}+(k+1)(2+\sqrt{3})\bigr]+\bigl(\sqrt{3}-2\bigr)L\sum_{j=1}^{k+1} t^{2j-2}\\
Lt^{2k+2}\bigl[\frac{1}{2\sqrt{3}}+1-(k+1) \bigr]+L\sum_{j=1}^{k+1} t^{2j-2}.
\end{array} \right.
\end{eqnarray}
Recalling that
$$
\sum_{j=1}^{\infty}t^{2j-2}=\frac{1}{1-t^2}
$$
it is now immediate to observe that $\widetilde{\uu}(x_{A_k},y_{A_k})$, $\widetilde{\uu}(x_{B_k},y_{B_k})$, $\widetilde{\uu}(x_{C_k},y_{C_k})$ and $\widetilde{\uu}(x_{A_{k+1}},y_{A_{k+1}})$ 
%
%
all converge to $\epsilon L(\sqrt{3}-2,1)\frac{1}{1-t^2}$ as $k\to\infty$. Repeating this argument for the remaining regions $\omega_{B_k},\dots, \omega_{F_k}$ yields
%
%
\begin{eqnarray}\label{1306252237}
\MM_k=(M_{1,k},M_{2,k})\to \epsilon L(\sqrt{3}-2,1)\frac{1}{1-t^2}
\quad \textrm{as} \quad k\to\infty.
\end{eqnarray}
The, we can define
\begin{eqnarray}\label{1306252225}
m_{i,k}:=\inf_{(x,y)\in\omega_{k}}\widetilde{u_i}(x,y)=\min\Bigl\{\Bigl(\inf_{(x,y)\in\omega_{A_k}}\widetilde{u_i}(x,y)\Bigr),\Bigl(\inf_{(x,y)\in\omega_{B_k}}\widetilde{u_i}(x,y)\Bigr),\nonumber\\\Bigl(\inf_{(x,y)\in\omega_{C_k}}\widetilde{u_i}(x,y)\Bigr), \Bigl(\inf_{(x,y)\in\omega_{D_k}}\widetilde{u_i}(x,y)\Bigr),
\Bigl(\inf_{(x,y)\in\omega_{E_k}}\widetilde{u_i}(x,y)\Bigr),\Bigl(\inf_{(x,y)\in\omega_{F_k}}\widetilde{u_i}(x,y)\Bigr)\Bigr\}.
\end{eqnarray}
Recalling that the infimum of $\widetilde{\uu}$ is attained on the vertices, the same computation of $(\ref{1306252230}\hbox{-}\ref{1306252234})$ is enough to prove that the sequence $\{\inf_{\omega_{A_k}}(\widetilde{u_1}),\inf_{\omega_{A_k}}(\widetilde{u_2})\}$ converges to $\epsilon L(\sqrt{3}-2,1)\frac{1}{1-t^2}$ as $k\to\infty$. Repeating the same argument for the remaining infima in (\ref{1306252225}) shows that
\begin{eqnarray}\label{1306252238}
\mm_k=(m_{1,k},m_{2,k})\to \epsilon L(\sqrt{3}-2,1)\frac{1}{1-t^2}
\quad \textrm{as} \quad k\to\infty.
\end{eqnarray}
Finally, by combining (\ref{1306252237}) and (\ref{1306252238}) we have 
\begin{eqnarray}
\widetilde{\uu}(x,y)\to \epsilon L(\sqrt{3}-2,1)\frac{1}{1-t^2}
\quad \textrm{as} \quad (x,y)\to(0,0)
\end{eqnarray}
as claimed.
Therefore, by a further extension of $\widetilde{\uu}$ we can define the continuous function
\begin{eqnarray}\label{uuuu}
\uu(x,y) :=\left\{
\begin{array}{cccc}
\widetilde{\uu}(x,y) &\textrm{on}\quad\Omega\slash (0,0)\\
\epsilon L(\sqrt{3}-2,1)\frac{1}{1-t^2} &\textrm{for } (x,y)=(0,0)\\
\textrm{extended by continuity on $\overline{\partial\Omega}$.}
\end{array} \right.
\end{eqnarray}
The fact that $\widetilde{\uu}_A(x_A,y_A)=\widetilde{\uu}_B(x_A,y_A)=\widetilde{\uu}_C(x_A,y_A)$, $\widetilde{\uu}_C(x_F,y_F)=\widetilde{\uu}_F(x_F,y_F)=\widetilde{\uu}_D(x_F,y_F)$, and 
$\widetilde{\uu}_D(x_E,y_E)=\widetilde{\uu}_E(x_E,y_E)=\widetilde{\uu}_B(x_E,y_E)$ guarantees that the function defined above is indeed continuous up to the boundary of $\Omega$.

 We now prove point 2). 
Since $\uu(x,y)$ is affine over each of the sets $\omega_{A_k}$, $\omega_{B_k}$, $\omega_{C_k}$, $\omega_{D_k}$, $\omega_{E_k}$ and $\omega_{F_k}$, it is straightforward to compute its gradient: 
\begin{eqnarray}\label{1307111146}
\nabla\uu(x,y) =\HH(x,y):=\left\{
\begin{array}{cccc}
\HH_{A_k}(x,y) &\textrm{on}\quad\omega_{A_k}\\
\HH_{B_k}(x,y) &\textrm{on}\quad\omega_{B_k}\\
\HH_{C_k}(x,y) &\textrm{on}\quad\omega_{C_k}\\
\HH_{D_k}(x,y) &\textrm{on}\quad\omega_{D_k}\\
\HH_{E_k}(x,y) &\textrm{on}\quad\omega_{E_k}\\
\HH_{F_k}(x,y) &\textrm{on}\quad\omega_{F_k}\\
\end{array} \right.
\end{eqnarray}
where
\begin{eqnarray}\label{1306262148}
\HH_{A_k}=
\epsilon\left(
\begin{array}{ccc}
-\frac{1}{2} & \sqrt{3}(2k-\frac{1}{2})  \\
-\sqrt{3}(2k+\frac{1}{2}) & \frac{1}{2}
\end{array} \right)=\EE_3+\WW_k
\end{eqnarray}
\begin{eqnarray}\label{}
\HH_{B_k}=
\epsilon\left(
\begin{array}{ccc}
1 & (-\sqrt{3}+2k\sqrt{3} )  \\
(\sqrt{3}-2k\sqrt{3}) & -1
\end{array} \right)=\EE_1+\WW_k+\widetilde{\WW}
\end{eqnarray}
\begin{eqnarray}\label{}
\HH_{C_k}=
\epsilon\left(
\begin{array}{ccc}
-\frac{1}{2} & \sqrt{3}(2k-\frac{1}{2})  \\
\frac{3\sqrt{3}}{2}-2k\sqrt{3} & \frac{1}{2}
\end{array} \right)=\EE_2+\WW_k+\widetilde{\WW}
\end{eqnarray}
\begin{eqnarray}\label{}
\HH_{D_k}=
\epsilon\left(
\begin{array}{ccc}
-\frac{1}{2} & 2k\sqrt{3}-\frac{3\sqrt{3}}{2} \\
(\frac{\sqrt{3}}{2}-2k\sqrt{3}) & \frac{1}{2}
\end{array} \right)=\EE_3+\WW_k+\widetilde{\WW}
\end{eqnarray}
\begin{eqnarray}\label{}
\HH_{E_k}=
\epsilon\left(
\begin{array}{ccc}
-\frac{1}{2} & (\frac{1}{2}+2k)\sqrt{3} \\
(\frac{1}{2}-2k)\sqrt{3} & \frac{1}{2}
\end{array} \right)=\EE_2+\WW_k
\end{eqnarray}
\begin{eqnarray}\label{1306262149}
\HH_{F_k}=
\epsilon\left(
\begin{array}{ccc}
1 & 2k\sqrt{3} \\
-2k\sqrt{3} & -1
\end{array} \right)=\EE_1+\WW_k
\end{eqnarray}
 and 
\begin{eqnarray}\label{}
\WW_{k}=
\epsilon\left(
\begin{array}{ccc}
0 & 2k\sqrt{3} \\
-2k\sqrt{3} & 0
\end{array} \right),\quad
\widetilde{\WW}=
\epsilon\left(
\begin{array}{ccc}
0 & -\sqrt{3}  \\
\sqrt{3} & 0
\end{array} \right).
\end{eqnarray}
Therefore we have $\bigl(\nabla\uu(x,y)\bigr)_{sym}\in \{\EE_1,\EE_2,\EE_3\}$ for almost every $(x,y)$ in $\Omega.$

Proof of point 3) is also straightforward. By combining (\ref{1307111146}) with (\ref{1306262148}-\ref{1306262149}) we can write
$$
\calL^2\Bigl((x,y)\in\Omega:\bigl(\nabla\uu(x,y)\bigr)_{sym}=\EE_1\Bigr)=\sum_{k=0}^{\infty}\Bigl(\calL^2(\omega_{F_k}+\calL^2(\omega_{B_k})\Bigr),
$$
$$
\calL^2\Bigl((x,y):\bigl(\nabla\uu(x,y)\in\Omega\bigr)_{sym}=\EE_2\Bigr)=\sum_{k=0}^{\infty}\Bigl(\calL^2(\omega_{E_k}+
\calL^2(\omega_{C_k})\Bigr),
$$
$$
\calL^2\Bigl((x,y):\bigl(\nabla\uu(x,y)\in\Omega\bigr)_{sym}=\EE_3\Bigr)=\sum_{k=0}^{\infty}\Bigl(\calL^2(\omega_{A_k}+\calL^2(\omega_{D_k})\Bigr).
$$
Then, recalling that by (\ref{1307111200}) we have $\calL^2(\omega_{A_k})=\calL^2(\omega_{E_k})=\calL^2(\omega_{F_k})$ and $\calL^2(\omega_{B_k})=\calL^2(\omega_{C_k})=\calL^2(\omega_{D_k})$ for every $k\geq 0$, we conclude that 
$$
\sum_{k=0}^{\infty}\Bigl(\calL^2(\omega_{F_k}+\calL^2(\omega_{B_k})\Bigr)=\sum_{k=0}^{\infty}\Bigl(\calL^2(\omega_{E_k}+
\calL^2(\omega_{C_k})\Bigr) =\sum_{k=0}^{\infty}\Bigl(\calL^2(\omega_{A_k}+\calL^2(\omega_{D_k})\Bigr)=\frac{\calL^2(\Omega)}{3}
=\frac{\pi}{3} R^2.
$$

We are left with the proof of point 4), that is, $\uu\in W^{1,p}(\Omega,\RRR^2)$ with $1\leq p<\infty$. We only have to show that $\nabla\uu\in L^{p}(\Omega,\MMM^{2\times 2})$  with $1\leq p<\infty$. Componentwise, we have to show that $\nabla u_{ij}\in L^{p}(\Omega)$  with $i,j=1,2$. 
For $(x,y)\in\omega_k$ we have
$$
|\nabla u_{ij}(x,y)|^p\leq 2^{p-1}\bigl(|\nabla u_{sym,ij}(x,y)|^p+|\nabla u_{skew,ij}(x,y)|^p\bigr).
$$
Notice that by (\ref{1306262148})-(\ref{1306262149}) $\nabla u_{sym,ij}\in L^{\infty}(\Omega)$ for $i,j=1,2$.
Therefore, to obtain the claim we only have to show that $\nabla u_{skew,ij}$ is in $L^p(\Omega)$ with $1\leq p<\infty$. Of course, it is enough to consider the case $i\neq j$. 
By Beppo Levi theorem on monotone convergence we have 
\begin{gather}
\int_{\Omega}|(\nabla u_{skew})_{ij}|^p dxdy=
\int_{\Omega}\sum_{k=0}^{\infty}|(\nabla u_{skew})_{ij}|^p\chi_{\omega_k}(x,y)\,dxdy=\nonumber\\
\lim_{N\to\infty}\sum_{k=0}^{N}\int_{\Omega}|(\nabla u_{skew})_{ij}|^p\chi_{\omega_k}(x,y)\,dxdy
=
\lim_{N\to\infty}\sum_{k=1}^{N}\int_{\Omega}|(\nabla u_{skew})_{ij}|^p\chi_{\omega_k}(x,y)\,dxdy+const.,\nonumber
\end{gather}
where the characteristic function $\chi_{\omega_k}$ is equal to $1$ on $\omega_k$ and $0$ in $\Omega\backslash\omega_k$. 
For $k\geq 1$, in view of (\ref{1306262148})\hbox{-}(\ref{1306262149}), we can write 
\begin{eqnarray}\label{1306301321}
c_1^2 k^2\chi_{\omega_k}(x,y)\leq |(\nabla u_{skew})_{ij}|^2\chi_{\omega_k}(x,y)
\leq c_2^2 k^2\chi_{\omega_k}(x,y).
\end{eqnarray}
Here $c_1$ and $c_2$
are two (positive) real numbers that do not depend on $k$. 
Therefore we have 
\begin{eqnarray}\label{1309101806}
\lim_{N\to\infty}\sum_{k=1}^{N}\int_{\Omega}|(\nabla u_{skew})_{ij}|^p\chi_{\omega_k}(x,y)\,dxdy\leq
\lim_{N\to\infty}\sum_{k=1}^{N}\int_{\Omega}
 c_2^p k^p \chi_{\omega_k}(x,y)\,dxdy.
\end{eqnarray}
To conclude it is enough to estimate the right hand side in (\ref{1309101806}). Recalling (\ref{1307081703}) we can write 
\begin{gather}\label{}
\lim_{N\to\infty} c_2^p\sum_{k=1}^{N}\int_{\Omega}
 k^p\chi_{\omega_k}(x,y)\,dxdy
=
\lim_{N\to\infty}3\,c_2^p\sum_{k=1}^{N}\Bigl(\int_{\omega_{A_k}}
 k^p \,dxdy+\int_{\omega_{B_k}}
 k^p \,dxdy\Bigr)=\nonumber\\
\lim_{N\to\infty}3\,c_2^pL^2\sum_{k=1}^{N}\Bigl( 
 t^{4k}t k^p + 
t^{4k}t^{-1} k^p\Bigr).
\end{gather}
Since $t=\tan(\pi/12)<1$ this series converges for $1\leq p<\infty$.
$\qquad\qquad\qquad\qquad\qquad\qquad\square$

\begin{Remark}\label{REMA}

Notice that from point 4) of Theorem \ref{THM} it automatically follows by the Sobolev imbedding \cite [Thm. 6.1-3]{Ciarlet} that $\uu\in C^{0,r}(\Omega,\RR^2)\cap C(\overline{\Omega},\RR^2)$, where $C^{0,r}(\Omega,\RR^2)$ is the space of H\"{o}lder continuous functions with $0<r<1$. However, notice that our proof of point 4) implicitly makes use of the fact that $\uu$ does not jump across the interfaces, which follows from 1). Moreover, our proof of 1) has the advantage of actually determining the value of the map $\uu$ in the origin.

Interestingly, from Point 4) of Theorem \ref{THM} we conclude that $\uu$ is not Lipschitz continuous in $\Omega$. To verify this, observe that by Eq. (\ref{1306301321}) we can not find an uniform bound for $\nabla\uu_{skew}$ in $\Omega$. 
This is a relevant difference with respect to the theory of simple laminates \cite{Ball87}, \cite{Muller}, \cite{BhattaBook} in which solutions of differential inclusions are usually sought in the space of piecewise affine amps.
However, given any positive $\delta$, $\uu$ is Lipschitz continuous over the smaller set $\Omega\backslash \overline{B_{\delta}(0,0)}$.

\end{Remark}

\section{Discussion}
\label{1307112003}



\begin{figure}[h!]
\centering%
\includegraphics[width=6.8cm]{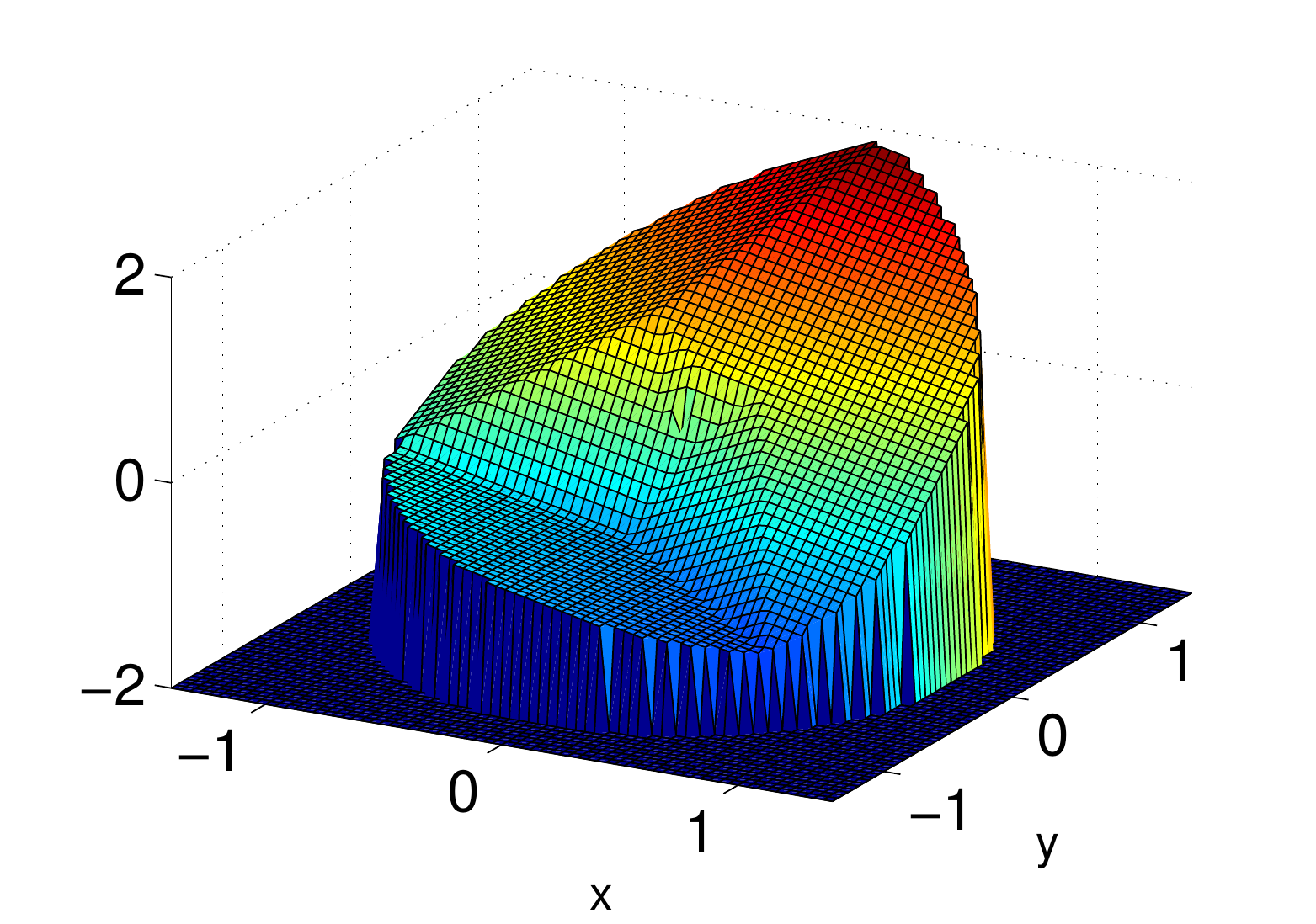}
\includegraphics[width=6.8cm]{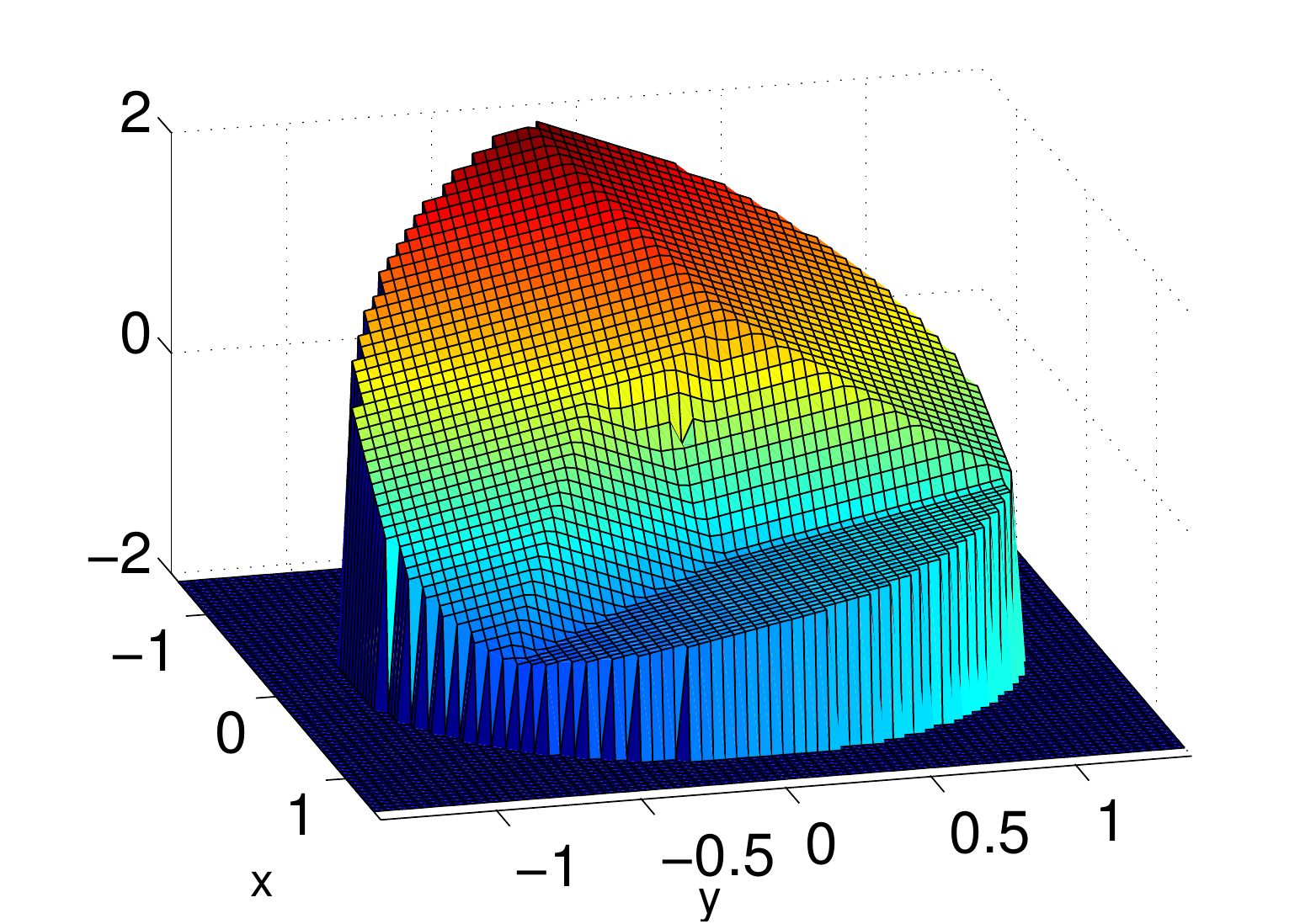}
\caption{Plot of $\uu(x,y)=(u_{\overline{z}1}(x,y),u_{\overline{z}2}(x,y))$. LEFT: $u_{\overline{z}1}(x,y)$. RIGHT:  $u_{\overline{z}2}(x,y)$ (here for simplicity L=$\epsilon$=1). The web version of this article contains the above plot figures in color.}\label{1011121630}
\end{figure}

\begin{figure}[h!]
\centering%
\includegraphics[trim = 50mm 00mm 50mm 00mm, clip=true, width=5.5cm ]{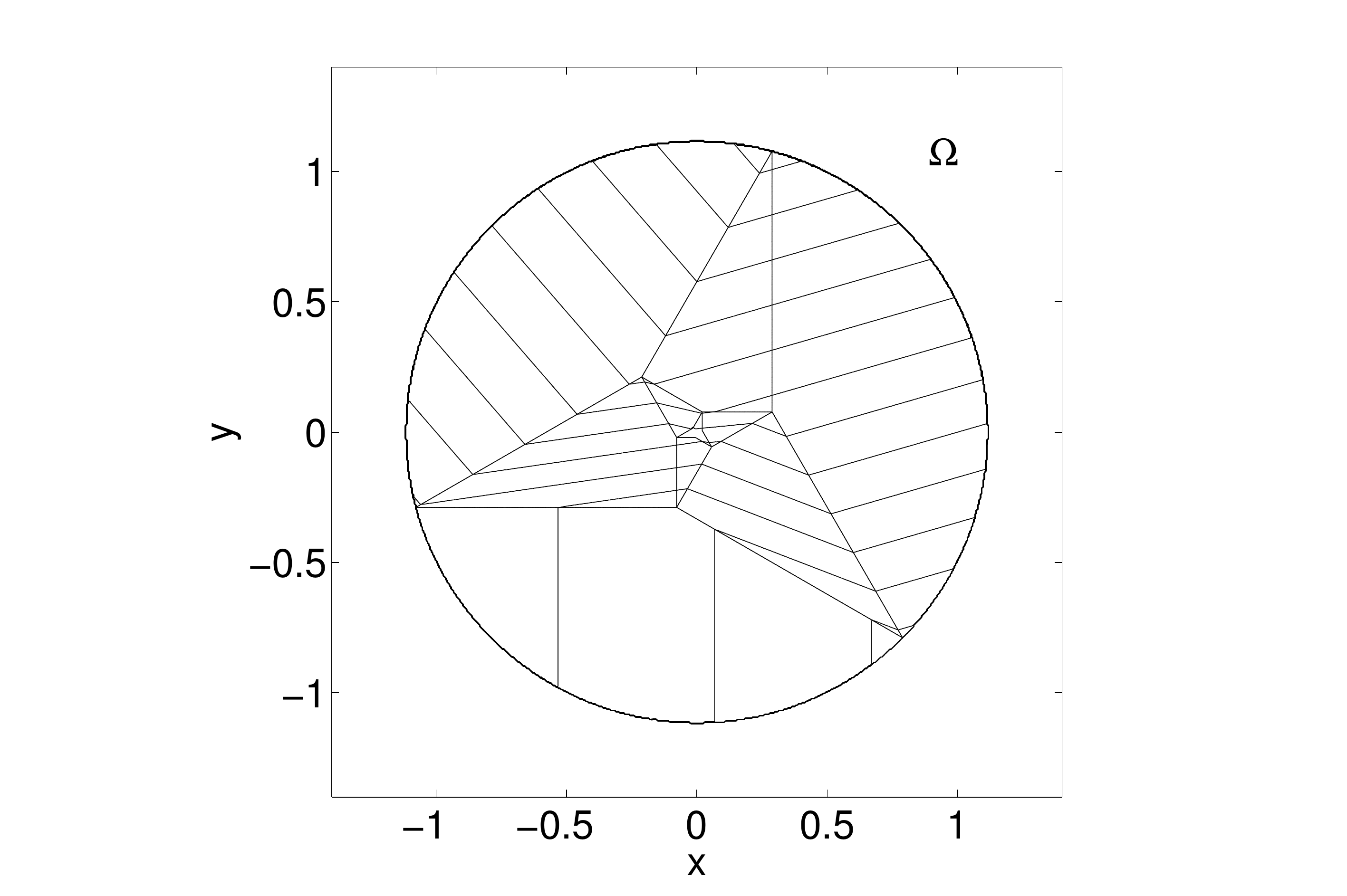}
\includegraphics[trim = 50mm 00mm 50mm 00mm, clip=true, width=5.5cm ]{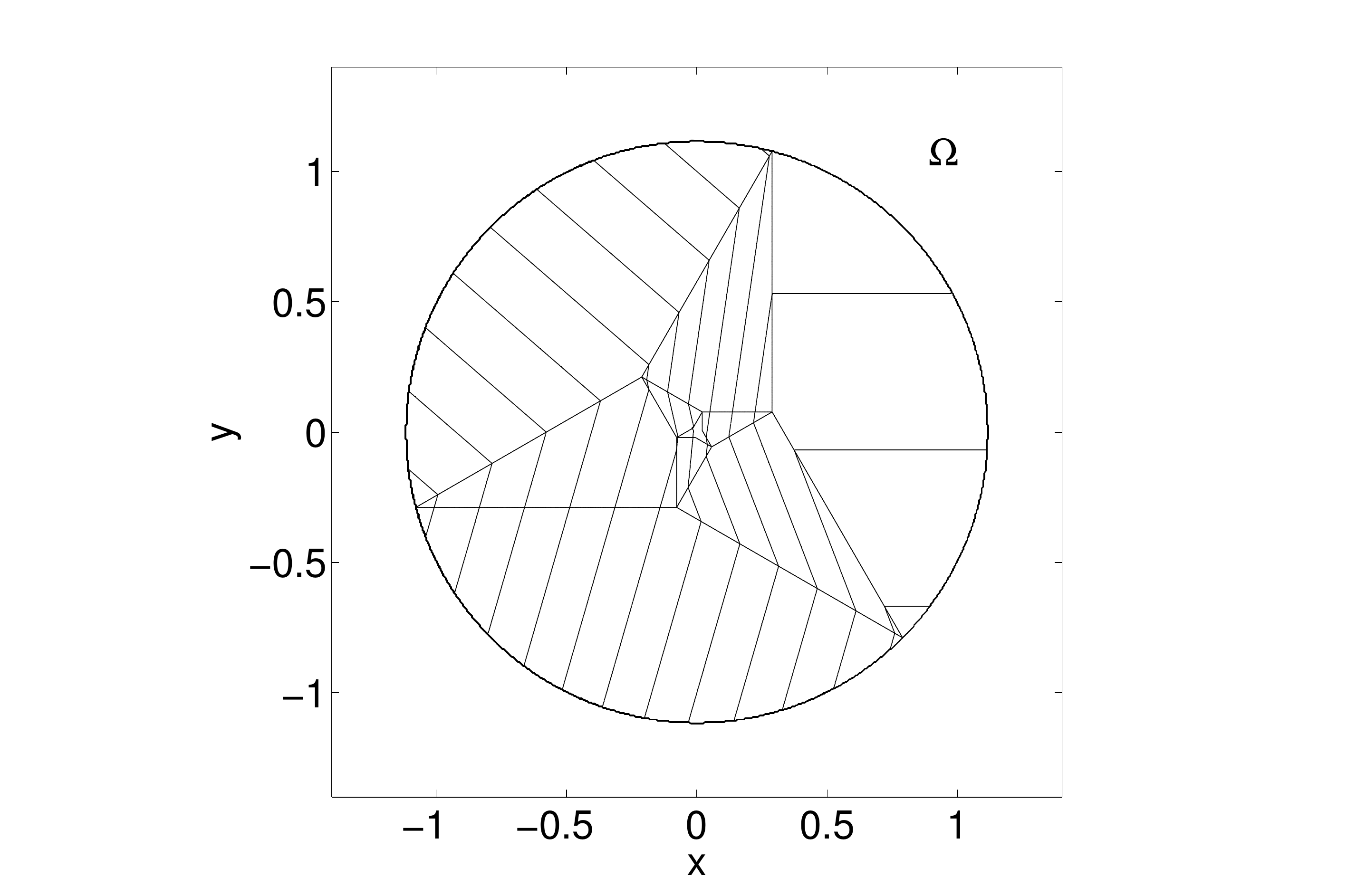}
\caption{Level curves of $\uu_{\overline{z}}(x,y)=(u_{\overline{z}1}(x,y),u_{\overline{z}2}(x,y))$. LEFT: $u_{\overline{z}1}(x,y)$. RIGHT: $u_{\overline{z}2}(x,y)$ (here for simplicity L=$\epsilon$=1). 
}\label{1011121631}
\end{figure}

Theorem \ref{THM} guarantees the existence of the function $\uu$ satisfying the properties stated in points 1-4) of the theorem, but it does not guarantee its uniqueness. 
The absence of additional constraints in the model, such as boundary conditions or applied forces, enables the construction of a different map $\uu$ still satisfying all the points of the theorem possibly with a different geometry. For instance,  the definition (\ref{uuuu}) can be modified by adding a rigid deformation (that is, a map whose gradient is a constant skew-symmetric matrix) and the new function $\uu$ still fulfills points 1-4) of Theorem \ref{THM}.
Indeed, let $\zz=(z_1,z_2,z_3)$ be any vector in $\RRR^3$ and let
\begin{eqnarray}\label{1306262257}
\uu_z(x,y)=\uu(x,y)+
\left(
\begin{array}{ccc}
0 & z_1 \\
-z_1 & 0
\end{array} \right) 
\left(
\begin{array}{ccc}
x \\
y
\end{array} \right)+
\left(
\begin{array}{ccc}
z_{2} \\
z_3 
\end{array} \right),
\end{eqnarray}
where $\uu$ is defined in (\ref{uuuu}). Then, $\uu_z(x,y)$ clearly fulfills points 1-4) of Theorem \ref{THM}.
Letting $z=\overline{z}=(3\frac{\sqrt{3}}{2}, 0, 0)$, in Fig. \ref{1011121630}  we plot the graph of the function $\uu_{\overline{z}}(x,y)$ and in Fig. \ref{1011121631} we plot the level curves of $\uu_{\overline{z}}$. The choice of $\overline{z}$ is purely for clear representation purposes.
From these figures it becomes clear that the field $\uu_{\overline{z}}(x,y)$ is locally piecewise affine (except at each neighborhood of the origin), and that there is no singularity at the center of the tripole-star pattern. In particular, it has to be noticed that the displacement $\uu_{\overline{z}}(x,y)$ is continuous across the interfaces. This property is a direct consequence of the geometrical compatibility condition which has been imposed in the construction of $\uu(x,y)$.
While the construction of Theorem $\ref{THM}$ works for any pair $L>0$, $\varepsilon>0$, it is clear that some limitations have to be imposed on $L,\varepsilon$ in order that the assumptions of linearized elasticity are fulfilled.
Since $\uu$ is continuous in $\overline{\Omega}$, it is straightforward to derive conditions on $L,\varepsilon$ that guarantee that $|\uu|=o(1)$. In fact, by denoting with $u_1(x,y)$ and $u_2(x,y)$ the components of $\uu(x,y)$,
we can estimate
$$
\max_{(x,y)\in\overline{\Omega}}|\uu(x,y)|\leq\sqrt{\Bigl[\max_{(x,y)\in\overline{\Omega}}|u_1(x,y)|\Bigr]^2+\Bigr[\max_{(x,y)\in\overline{\Omega}}|u_2(x,y)|\Bigr]^2}
$$
where
$$
\max_{(x,y)\in\overline{\Omega}}|u_i(x,y)|=\max\Bigl\{\max_{(x,y)\in\overline{\Omega}}u_i(x,y),\,\,-\min_{(x,y)\in\overline{\Omega}}u_i(x,y)\Bigr\},\,\,\, i=1,2.
$$
As in the proof of Point 1) of Theorem \ref{THM}, we are left with splitting the computation of the maximum and the minimum of the components of $\uu$ over $\overline{\Omega}$ into the evaluation of $u_1(x,y)$ and $u_2(x,y)$ on the vertices of the rhomboids:
\begin{eqnarray}\label{1307151513}
\max_{(x,y)\in\overline{\Omega}}u_i(x,y)=\max\Bigl\{u_i(x_{A_k},y_{A_k}),u_i(x_{B_k},y_{B_k}), u_i(x_{C_k},y_{C_k}),u_i(x_{D_k},y_{D_k}),\qquad\qquad\quad\nonumber\\
\qquad\qquad\qquad\quad u_i(x_{E_k},y_{E_k}),u_i(x_{F_k},y_{F_k}), u_i(0,0),\quad k\in \{0\}\cup\NNN\Bigr\},\,\,\, i=1,2,
\end{eqnarray}
\begin{eqnarray}\label{1307151514}
\min_{(x,y)\in\overline{\Omega}}u_i(x,y)=\min\Bigl\{u_i(x_{A_k},y_{A_k}),u_i(x_{B_k},y_{B_k}), u_i(x_{C_k},y_{C_k}),u_i(x_{D_k},y_{D_k}),\qquad\qquad\quad\nonumber\\
\quad\qquad\qquad\quad u_i(x_{E_k},y_{E_k}),u_i(x_{F_k},y_{F_k}), u_i(0,0),\quad k\in \{0\}\cup\NNN\Bigr\},\,\,\, i=1,2.
\end{eqnarray}
Now, by simply plugging the equations of the vertices (\ref{1307151522}) into (\ref{uuuu}) we immediately obtain that
\begin{eqnarray}\label{1307151513}
\max_{(x,y)\in\overline{\Omega}}u_i(x,y)\propto\epsilon L.
\end{eqnarray}
Therefore, if $\epsilon L \ll 1$ we have that  $|\uu(x,y)|=o(1)$ for every $(x,y)$ in $\overline{\Omega}$.

The discussion for the gradient of the displacement is more delicate.
The fact that $\nabla\uu_{skew}$ is unbounded in $\Omega$ (see Remark \ref{REMA}) has deep consequences on the validity range of our model in the scenario of linearized elasticity (see \cite{BhattaBook} Chapter 11). Correctly, such a model holds under the assumption that $|\nabla\uu|$ is small. This poses a requirement on the symmetric  and skew symmetric parts of $\nabla\uu$ that have to be simultaneously small. While by Eq. (\ref{1306261933})
 $|\nabla\uu_{sym}(x,y)|=\sqrt{2}\epsilon$ is uniformly bounded at almost every $(x,y)$ in $\Omega$, Eq. (\ref{1306301321}) tells us that $|\nabla\uu_{skew}|$ grows linearly in $k$ over $\omega_k$. 
Therefore for $k$ large enough, the assumptions of geometrically linearized elasticity fail.
Whether replacing the linearized theory with the full non-linear elasticity model would bring more insight onto the analysis close to the origin is not clear.
Although the full analysis of nested microstructure of the type of Fig.\ref{1305141104} in the more general model of finite elasticity is to date an open question, here we try to address this point with an example. 
Consider for instance the point 
$$\xi_k:=(x_{\xi_k},y_{\xi_k})=\frac{1}{2}\bigl(x_{B_k}+x_{B_{k+1}},y_{B_k}+y_{B_{k+1}}\bigr)=Lt^{2k}\Bigl(\frac{5}{\sqrt{3}}-3, 3-\frac{5}{\sqrt{3}}\Bigr)$$ 
which belongs to $\omega_{B_k}$ and is equally distant from $B_k$ and $B_{k+1}$. As $k\to\infty$ the sequence $\{\xi_k\}$ converges to $(0,0)$ at a rate of $t^{2k}$.
%
%
%
%
In the case where $L\approx 0.07 \mu$m (as can be
deduced from Fig. \ref{1305141104}-(a)) we have that, already for $k=2$
$$
|\xi_2|\sim 10^{-2} nm \hspace{2cm} 
$$
which is a distance smaller than the crystal lattice parameters. 
Thus, in order to capture all the details of the behavior of the system close to the origin, even the use of a full non-linear elastic model
would not be justified. This is where  atomistic models would be more reliable.
Nevertheless, the domains obtained for $k=0,1$ and corresponding to the regions $\omega_0=\omega_A\cup\omega_B\cup\omega_C\cup\omega_D\cup\omega_E\cup\omega_F$
and  $\omega_{B_1},\omega_{C_1}$ and $\omega_{D_1}$ 
remain in the range for which continuum approximations are acceptable. Notice that due to the fast convergence in $k$, these regions cover the great majority of our reference domain $\Omega$. Indeed, a straightforward computation shows that the area of this region corresponds to
$$
\frac{\calL^2(\omega_0)+3\calL^2(\omega_{B_1})}{\calL^2(\Omega)}\approx 0.999,
$$
thus confirming that the region 
corresponding to $k\geq 2$
is confined in a tiny neighborhood of the origin.

We want to remark that the domain-length $L=0.07 \mu$m reported above
represents only a possible example of an observed microstructure. Due to self-similarity, other tripole-stars with larger $L$ may well be observed for which higher $k$-order generations remain in the range where continuum models are valid.

Finally, we discuss the effect of neglecting the interfacial energy. 
Computations based on a model for the triangle-to-centered rectangle transformation with sharp interfaces (A. Ruland, forthcoming) show that minimizers of the total energy may develop singular gradients also in the presence of interfacial energy terms. 
Formation of the self-similar microstructure is not prevented by an energy penalisation of the interfaces because the line energy contribution associated with the nesting is negligible.
This energy argument 
thus validates our findings based on a model with no interface energy contribution.

\quad
\\
\noindent In summary,  we have analyzed a nested martensitic microstructure observed in lead orthovanadate and Mg-Cd alloys.
In the scenario of linearized elasticity, we have modeled the microstructure as a three-phase martensitic mixture by showing that geometric compatibility across the interfaces holds for each interface between the variants making up the hierarchical microstructure.
Continuity of the displacement (Theorem \ref{THM}, point 1)) rules out irreversible phenomena such us cracks and cavitations in the model analyzed in this article. 
This work also offers an analytical counterpart to the numerical work of Ref. \cite{Lookman01}. 
The approach adopted in Ref.  \cite{Lookman01} is \textit{strain-based}, in that the minimization problem is written in terms of strain components (rather than on the displacement map $\uu$) and geometric compatibility is imposed through a partial differential equation relating the strain components.
Thus,  the properties  of $\uu$ are not investigated. 
In Theorem 1, we have proved that  geometric compatibility is not only achieved globally in the domain but also close to the origin. We have also shown 
 that the map $\uu$, which realizes the microstructure, is smooth.

 One of the main contributions of this paper is the precise estimation for the growth rate of the gradient of the displacement and on the computation of the geometry of the domains
close to the singularity. 
The unboundedness of the gradient of the displacement is an intrinsic aspect of the analyzed system.
While no continuum model, either in non-linear or geometrically linearized elasticity, is appropriate to exactly describe the singularity, we provide precise quantitative information on the scaling law of the self-similar structure and estimate the size of the region where a linearized model becomes inconsistent.
Remarkably, the geometry of the microstructure we compute corresponds to the one observed experimentally and to the ones computed numerically with a good level of approximation (see Fig. \ref{1305141104}).

\newpage
Whether the results of this paper can be extended to treat more general situations belonging to  the large set of non-periodic microstructures and, possibly fractals, is an intriguing question. This will be considered in a forthcoming paper. However, we have addressed one of the most challenging mathematical aspects associated with this question. 
Nesting of microstructure close to a singleton may turn into an increasing gradient of the displacement.
%
%
%
%
Although in the current situation this results in loss of regularity of $\uu$ which fails to be Lipschitz in $\Omega$), it is not clear if 
the approach described in this paper may be applied to more general situations. 
%
%
Thus, a more general theory of nested martensitic microstructure requires treatment of geometric compatibility in the presence of potential point defects. Domains with singularities may be seen as  special cases of non-simply-connected sets and therefore would require the types of techniques developed in Ref. \cite{Yavari} (see also \cite[Chapter 2]{Barber}).

\subsection*{Acknowledgments}
We acknowledge support from the Department of Energy
National Nuclear Security Administration under Award Number
DE-FC52-08NA28613.
P.C. is grateful to Los Alamos National
Laboratory for its kind hospitality.
P.C. was partially supported by the European Research Council under the European Union's Seventh Framework Programme (FP7/2007-2013) - ERC grant agreement N. 291053.
Part of this work has been written when P.C. was a postdoctoral student at
California Institute of Technology.
The authors are grateful with Kaushik Bhattacharya, Richard James and Angkana Ruland for several discussions at various stages of the work.


\begin{thebibliography}{00}

\bibitem{Ball87} J.M. Ball, R.D. James, 1987. Fine phase mixtures as minimizers of energy. Arch. Rat. Mech. Anal., 100, 13-52.
\bibitem{Barber} J.R. Barber, 2002. Elasticity. Kluwer, Dordrecht.
\bibitem{BhattaBook} K. Bhattacharya, 2003. Microstructure of martensite. Oxford University Press.
\bibitem{Ciarlet} P.G. Ciarlet, 1988. Mathematical
elasticity Vol.1. Elsevier Science Publishers B.V.
\bibitem{curnoe} S. H. Curnoe,  A. E. Jacobs, 2001. Phys Rev B 63, 094110.
\bibitem{Jacobs} A.E. Jacobs, S.H. Curnoe, R.C. Desai, 2004. Landau Theory of domain patterns in ferroelastic. Mater. Trans. 45, 1054-1059.
\bibitem{Khachachuryan} A. G. Khachaturyan, 1983. Theory of structural transformations in solids. John Wiley and Sons, NY.
\bibitem{kitano} Y. Kitano, K. Kifune, 1991. HREM study of disclinations in MgCd ordered alloy. Ultramicroscopy 39, 279-286.
\bibitem{Lookman03} T. Lookman, S.R. Shenoy, K. O. Rasmussen, A. Saxena, A.R. Bishop, 2003. Ferroelastic dynamics and strain compatibility. Phys. Rev. B 67.
\bibitem{Manolikas} C. Manolikas, S. Amelinckx, 1980. Phase transitions in ferroelastic lead orthovanadate as observed by means of electron microscopy and electron diffraction. Phys. Stat. Sol. (a) 60, 607-617.
\bibitem{Muller} S. M\"{u}ller, 1999. Variational methods for microstructure and phase transitions, in: Proc. C.I.M.E. summer school ``Calculus of variation and geometric evolution problems'', Cetraro 1996, (F. Bethuel, G. Huisken, S. M\"{u}ller, K. Steffen, S. Hildebrandt, M. Str\"{u}we eds.),
Springer LNM vol. 1713.
\bibitem{Lookman01} M. Porta, T. Lookman, 2013. Heterogeneity and phase transformation in materials: energy minimization, iterative methods and geometric nonlinearity. Acta Materialia, Volume 61, Issue 14, Pages 5311-5340.  
\bibitem{Yavari} A. Yavari, 2013. Compatibility equations of nonlinear elasticity for non-simply connected bodies. Arch. Rat. Mech. Anal., 209, 237-253.
\end{thebibliography}
\end{document}